\documentclass[journal]{IEEEtran}
\usepackage{cite}
\usepackage{balance}
\usepackage{graphicx}
\usepackage{amsmath}
\usepackage{amssymb}
\newcommand{\comment}[1]{}

\newcommand{\halmos}{\rule{1ex}{1.4ex}}
\newcommand{\qed}{\hfill \mbox{$\halmos$}}

\newcommand{\kl}{{\mathcal KL}}
\newcommand{\kk}{{\mathcal K}}
\newcommand{\ki}{{\mathcal K}_\infty}

\newcommand{\mer}{\hfill $\circ$}

\def \R{{\mathbb R}}
\def \N{{\mathbb N}}

\def\I{\mathcal{I}}
\def\X{\mathcal{X}}

\def\K{\mathcal{K}}
\def\kl{\mathcal{KL}}
\def\S{\mathcal{S}}

\def\T{\mathcal{T}}

\def\U{\mathcal{U}}

\newtheorem{teo}{Theorem}[section]
\newtheorem{itlemma}{Lemma}[section]
\newtheorem{lema}{Lemma}[section]
\newtheorem{corol}{Corolary}[section]
\newtheorem{defin}{Definition}[section]

\newtheorem{itremark}[itlemma]{Remark}
\newtheorem{as}{Assumption}
\newtheorem{claim}{Claim}
\newenvironment{remark}{\noindent\begin{itremark}\rm}{\end{itremark}}
\newtheorem{ex}{Example}[section]

\begin{document}

\title{Uniform asymptotic stability of switched nonlinear time-varying  systems and detectability of reduced limiting control systems}
\author{J. L.~Mancilla-Aguilar 
        and~R. A. Garc\'{\i}a% <-this % stops a space
\thanks{J. L. Mancilla-Aguilar R. A. Garc\'{\i}a are with the Center of Systems and Control, Instituto Tecnol\'ogico de Buenos Aires, Buenos Aires, Argentina. (e-mail: \texttt{jmancill@itba.edu.ar, 
		ragarcia@itba.edu.ar})}

 \thanks{Work partially supported by ANPCyT grant PICT 2014-2599, Argentina.}}

\maketitle

\begin{abstract}
This paper is concerned with the study of both, local and global, uniform asymptotic stability for switched nonlinear time-varying (NLTV) systems through the detectability of output-maps. With this aim the notion of reduced limiting control systems for switched NLTV systems whose switchings verify time/state dependent constraints, and the concept of weakly zero-state detectability for those reduced limiting systems are introduced. Necessary and sufficient conditions for the (global)uniform asymptotic stability of families of trajectories of the switched system are obtained in terms of this detectability property. These sufficient conditions in conjunction with the existence of multiple weak Lyapunov functions, yield a criterion for the (global) uniform asymptotic stability of families of trajectories of the switched system. This criterion can be seen as an extension of the classical Krasovskii-LaSalle theorem. An interesting feature of the results is that no dwell-time assumptions are made. Moreover, they can be used for establishing the global uniform asymptotic stability of switched NLTV system under arbitrary switchings. The effectiveness of the proposed results is illustrated by means of various interesting examples, including the stability analysis of a semi-quasi-Z-source inverter.  
\end{abstract}

\begin{IEEEkeywords}
Switched nonlinear time-varying systems, asymptotic stability, weak zero-state detectability, limiting systems
\end{IEEEkeywords}

\IEEEpeerreviewmaketitle
 \section{Introduction}
The study of switched systems is a topic of current interest, mainly due to their versatility in the modeling of many instances of complex engineering systems  \cite{Liberzonbook,Liberzon-Morse, matveev, van-der-Schaft}. This versatility nevertheless comes at a price: the study of switched systems, in particular of their stability, may be rather involved. In fact, the behavior of switched systems may become very rich --and consequently difficult to characterize--, due to the interaction between the switching signals and the dynamics of the component subsystems. This fact motivated an  extensive investigation of different stability properties of switched systems (see \cite{Liberzonbook, Liberzon-Morse, DeCarlo, Shorten, Lin} and references therein).

In this regard, in  this work we intend to demonstrate the (global) uniform asymptotic stability, (G)UAS for short, of a large class of switched nonlinear time-varying (NLTV) systems allowing either arbitrary switching or state and/or time dependent switching.

It is well-known that the existence of common or multiple %strict 
Lyapunov functions suffices to state the (G)UAS of switched systems under, respectively, arbitrary switchings or restricted ones \cite{Liberzonbook}.  
Nevertheless, their obtention even for nonswitched systems is in general a highly non-trivial task. On the other hand, common or multiple weak Lyapunov functions, i.e. functions for which the time derivatives along the trajectories of the subsystems are negative semidefinite, arise rather naturally in many engineering systems. In fact, in many cases, the subsystems have dissipative models and the energy functions usually constitute common or multiple weak Lyapunov functions. 
Although the existence of weak Lyapunov functions by itself gives no asymptotic stability guarantee, it can be supplemented with other additional conditions in order to yield asymptotic stability. This fact motivated the development of several stability results for switched time-invariant systems based on weak Lyapunov functions: extensions of LaSalle's invariance principle 
\cite{hespan_tac04,bacmaz_scl05,chewan_tac08,goesan_scl08,wanche_scl09, zhajia_tac14, mancilla-garcia-3, mancilla-garcia-2} and other approaches  \cite{serres, riedinger, balde, jouan}. Most of these papers assume that the switching signals satisfy some dwell-time condition and only provide convergence results.   

Opposite to the time-invariant case, results based on weak Lyapunov functions valid for switched time-varying systems are scarce. A GUAS result was obtained in \cite{manhai_tac17} by using a perturbation approach. A refinement of Matrosov's theorem can be found in \cite{teel_matrosov_automatica16}. 
In the very interesting paper \cite{leejia_tac08}, extensions of the classical Krasovskii-LaSalle theorem were obtained, without assuming any dwell-time requirement,  by introducing the output persistently exciting (OPE) condition, whose verification is one of the main difficulties of this approach.

This consequently motivated the search of easier to check sufficient conditions for the OPE one (\cite{lee_tan_nesic_cdc2011, lee_tan_nesic_tac15, lee_tan_mareels_ecc16, lee_tan_mareels_tac17, lee_tan_mareels_cdc17}). In this regard, the concept of weak zero-state detectability (WZSD) was introduced in \cite{lee_tan_mareels_ecc16} and applied to switched NLTV systems. Although easier to check than the OPE condition, to determine the WZSD is still challenging. With the aim of simplifying the task of verifying WZSD, in \cite{lee_tan_mareels_tac17} the common zeroing-output system (CZOS) concept was introduced. Nevertheless, a CZOS does not always exist and even when it exists, to verify the WZSD is by no means trivial since the output signal depends on a sequence of switching signals.   

This paper is devoted to the study of the WZSD of switched NLTV systems with switched time-varying output and with switching signals satisfying some type of time/state dependent constraint, with the aim of obtaining simplified sufficient conditions for that detectability property (and consequently for the OPE condition). For doing that, and inspired by the notions of limiting solutions (see \cite{Artstein76} and \cite{Artstein78}) and of reduced limiting systems (see \cite{leejia_tac05}) for nonswitched NLTV systems, we introduce the concepts of zeroing-output limiting trajectories and of limiting reduced control systems. Roughly speaking, these zeroing-output limiting trajectories are the limits, in a suitable sense, of sequences of admissible trajectories of the switched system for which the output converges to zero, while the reduced limiting control systems are auxiliary control systems with outputs constrained to be zero and with inputs constrained to satisfy certain restrictions induced by those satisfied by the switching signals of the switched system. The key fact that {\em the zeroing-output limiting trajectories are solutions of these reduced limiting control systems} allows us to prove that the WZSD of the reduced limiting control systems, defined in a suitable way, implies that of the original switched system. Since in general it is easier to determine the WZSD of a reduced limiting control system than that of a switched system, we obtain in this way a simplified criterion for the WZSD of switched NLTV systems. 
This WZSD criterion in conjunction with other conditions, for example the existence of weak Lyapunov functions, allow us to yield new both sufficient and necessary conditions for the (G)UAS of switched NLTV systems. An interesting feature of these conditions is that they do not involve any dwell-time assumption.

The contribution of the paper can be summarized as follows:\\
\textbullet$\;$ New concepts of  reduced limiting control systems and of WZSD for these systems are introduced. These systems incorporate the time/state dependent constraints that the switching signals satisfy and their WZSD is far easier to check than that of the original switched system and implies the WZSD of the latter (Theorem \ref{thm:criterio_wzsd}). \\
\textbullet$\;$ A new characterization of the GUAS of a family of trajectories of the switched system in terms of the WZSD of its reduced limiting control systems is given (Theorem \ref{thm:converse}).\\
\textbullet$\;$ A criterion for (G)UAS of a family of trajectories, which can be seen as a generalization of the Krasovsii-LaSalle theorem, is obtained assuming the existence of weak multiple Lyapunov functions (Theorem \ref{thm:guas1}). It is also given a criterion, based on the existence of a common weak Lyapunov function,  for the GUAS of the trajectories of switched system under arbitrary switchings (Corollary \ref{cor:guas_as}).\\
\textbullet$\;$ The results that we present are by no means generalizations of those already obtained for switched nonlinear time-invariant systems, and since their reduced limiting control systems always exist and are easily obtained, we provide new tools for the stability analysis of these systems. 

The paper is organized as follows. In Section \ref{bdp} we state the problem addressed and introduce some of the concepts and assumptions employed. In Section \ref{sec:wzsd} we introduce the concepts of zeroing-output limiting trajectories, reduced limiting control systems, WZSD for these systems, and prove that this WZSD implies that of the original switched system. Section \ref{sec:guas} contains a characterization of the GUAS of a family of trajectories in terms of the WZSD of its limiting control systems. We present criteria for the UAS and the GUAS of the switched system under switching restrictions and arbitrary switching in Section \ref{sec:criteria}. Section \ref{ej} includes three interesting examples that illustrate the effectiveness of our results, while Section \ref{sec:conclusion} contains some conclusions. Finally, the Appendix provides the proofs of some auxiliary results.

\textbf{Notation:}  $\N$, $\R$ and $\R_{\ge 0}$ denote the sets of natural, real and nonnegative real numbers, respectively. For $x\in \R^n$, $|x|$ denotes its Euclidean norm. An indexed family of sets $\chi=\{\chi_i\}_{i\in \I}$ is a closed covering of $\R^n$ if $\chi_i$ is a closed subset of $\R^n$ for each $i \in \I$ and $\R^n=\cup_{i\in \I}\chi_i$.  Given $A\subset \R^m$, ${\rm co}(A)$, $\overline{A}$, ${\rm int}(A)$, $\partial A$ and $|A|$ are, respectively, the convex hull, the closure, the interior, the boundary and the Lebesgue measure of $A$. For any matrix $B$, $B'$ denotes its transpose. For any interval $J\subset \R$ and any $N \in \N$, $L^1_N(J)$ ($L^{\infty}_N(J)$) is the set of Lebesgue measurable functions $f:J\to \R^N$ which are Lebesgue integrable (essentially bounded). A function $\varphi:\Omega\to \R$, with $\Omega\subset \R^k$ for some $k$ and $0\in \Omega$, is said positive definite if $\varphi(0)=0$ and $\varphi(\omega)>0$ for all $\omega\neq0$. For a continuous function $\alpha:\R_{\ge 0} \to \R_{\ge 0}$, we write $\alpha\in\K$ if it is positive definite and strictly increasing, and $\alpha\in\ki$ if, in addition, $\alpha$ is unbounded. We write $\beta\in\mathcal{KL}$ if $\beta:\R_{\ge 0}\times \R_{\ge 0}\to \R_{\ge 0}$, $\beta(\cdot,t)\in\ki$ for any $t\ge 0$ and, for any fixed $r\ge 0$, $\beta(r,t)$ monotonically decreases to zero as $t\to \infty$.   
A function $g:J\times \X \to \R^m$, $(t,\xi) \mapsto g(t,\xi)$, with $J$ a real interval and $\X \subset \R^n$ is: i) uniformly bounded if there exists $J'\subset J$ with $|J\setminus J'|=0$ such that $g$ is bounded on $J'\times K$ for each compact subset $K\subset \X$; ii) measurable in $t$ if $g(\cdot,\xi)$ is Lebesgue measurable for all $\xi \in \X$; iii) continuous in $\xi$ if $g(t,\cdot)$ is continuous for almost all $t\in J$; iv) continuous in $\xi$ uniformly in $t$ if there exists a set $J'$ as in i) so that for every compact set $K\subset \X$ there exists a modulus of continuity $\omega_K\in \kk$ such that
$$ |g(t,\xi)-g(t,\zeta)|\le \omega_K(|\xi-\zeta|),\quad \forall t \in J',\;\forall \xi,\zeta \in K;$$  
v) locally Lipschitz in $\xi$ uniformly in $t$ if the modulus of continuity $\omega_K(r)\equiv L_K r$ for some $L_K\ge 0$.

\section{Problem statement}
\label{bdp}
For the index set $\I=\{1,\ldots,N\}$ and the open neighborhood of the origin $\X\subset \R^n$, we consider the switched NLTV system with outputs
\begin{eqnarray}
\dot{x}(t)&=&f(t,x(t),\sigma(t)) \label{eq:ss}\\
y(t)&=&h(t,x(t),\sigma(t)) \label{eq:so}
\end{eqnarray}
where, for $t\ge 0$, $x(t)\in \X$, $y(t)\in \R^p$ and $\sigma: \R  \rightarrow \I$ is a {\em switching signal}, {\em
i.e.}  $\sigma$ is piecewise constant (it has at most a finite number of jumps in each compact interval) and is continuous from the right. 

It is assumed throughout the paper that the system function $f:\R_{\ge 0}\times \X \times \I \to \R^n$ and the output-map $h:\R_{\ge 0} \times \X \times \I \to \R^p$ are such that for every $i\in \I$, $f_i(t,\xi):=f(t,\xi,i)$ and $h_i(t,\xi):=h(t,\xi,i)$ are measurable in $t$ and continuous in $\xi$, and $f_i$ satisfies the Carath\'eodory conditions \cite[p.28]{hale_book80}. 

Let $\S$ be the set of all the switching signals. A (forward) solution of (\ref{eq:ss}) corresponding to $\sigma \in \S$ is a locally absolutely continuous function $x:[t_0(x),t_f(x))\to \R^n$, with $0\le t_0(x)<t_f(x)$, such that $\dot{x}(t)=f(t, x(t),\sigma(t))$ for almost all $t\in [t_0(x),t_f(x))$. Due to the assumptions we made on $f$, for each $t_0\ge 0$, $x_0\in \X$ and $\sigma \in \S$ there is a (not necessarily unique) solution $x:[t_0(x),t_f(x))\to \X$ of (\ref{eq:ss}) such that $t_0(x)=t_0$, $t_f(x)>t_0(x)$ and $x(t_0)=x_0$. Such a solution is maximal if there does not exist a solution $\tilde x:[t_0(\tilde x),t_f(\tilde x))\to \R^n$ of (\ref{eq:ss}) corresponding to $\sigma$ such that $t_0(x)=t_0(\tilde x)$, $t_f(\tilde x)>t_f(x)$ and $\tilde x(t)=x(t)$ for all $t\in [t_0(x),t_f(x))$. We say that a pair $(x,\sigma)$ is a trajectory of (\ref{eq:ss}) if $x$ is a maximal solution of (\ref{eq:ss}) corresponding to $\sigma \in \S$. The trajectory $(x,\sigma)$ is forward complete if $t_f(x)=\infty$. 

The main concern of this paper is the study of the (global) uniform asymptotic stability of families of trajectories of system (\ref{eq:ss}). For completeness, we provide precise definitions of the stability properties we consider.
\begin{defin}\label{def:stab} \rm Let $\T$ be a family of trajectories of (\ref{eq:ss}).
 \begin{enumerate}
  \item $\T$ is uniformly stable (US) if there exist $\alpha \in \kk$ and $\delta_0>0$ such that
  for any $(x,\sigma)\in \T$ and any $s\in [t_0(x),t_f(x))$ such that 
  $|x(s)| \le  \delta_0$, we have $|x(t)|\le \alpha(|x(s)|)$ for all $s \le t <t_f(x)$.
  \item $\T$ is globally uniformly stable (GUS) if $\X = \R^n$ and there exists $\alpha \in \ki$ such that
  for any $(x,\sigma)\in \T$ and any $s\in  [t_0(x),t_f(x))$, we have $|x(t)|\le \alpha(|x(s)|)$ for all $s \le t <t_f(x)$.
  \item $\T$ is uniformly asymptotically stable (UAS) if there exist $\beta \in \mathcal{KL}$ and $\delta_0>0$ such that 
  for any $(x,\sigma)\in \T$ and any $s\in  [t_0(x),t_f(x))$ such that $|x(s)|< \delta_0$, we have $|x(t)|\le \beta(|x(s)|,t-s)$ for all $s \le t <t_f(x)$.
   \item \label{def:guas} $\T$ is globally uniformly asymptotically stable (GUAS) if $\X = \R^n$ and there exists $\beta \in \mathcal{KL}$ such that for any $(x,\sigma)\in \T$ and any $s\in  [t_0(x),t_f(x))$, we have $|x(t)|\le \beta(|x(s)|,t-s)$ for all $s \le t <t_f(x)$.
 \end{enumerate}
 \end{defin}
 \begin{remark}\label{rm:equiv}
  The stability properties defined above are equivalent to those defined in \cite[Defn.1]{leejia_tac08} in the classical $\varepsilon-\delta$ form. It is clear that the US property implies the forward completeness of any trajectory of $\T$ whose initial condition belongs to some open neighborhood $\mathcal{O}$ of the origin, and that the GUS property implies that the trajectories in $\T$ are forward complete. \mer
  \end{remark}

The US of $\T$ can be straightforwardly established if it admits common or multiple weak Lyapunov functions
\cite[Prop.1]{leejia_tac08}. Once the US of $\T$ is checked, the UAS of $\T$ can be obtained under additional conditions. A set of such additional conditions involves the use of an auxiliary output --which is often related to the total derivative of the weak Lyapunov function used for establishing the US of $\T$-- and the concept of weak zero-state detectability which we next recall.
  
We first introduce the notion of $\T$-zeroing-output sequence.
\begin{defin} \rm \label{def:zeroing-sequence} Given a family $\T$ of trajectories of (\ref{eq:ss}), the sequence $\{(t_k,x_k,\sigma_k)\}$ is a $\T$-zeroing-output sequence for (\ref{eq:ss})-(\ref{eq:so}) if $t_k\to \infty$, $(x_k,\sigma_k)\in \T$ for all $k$ and:
	\begin{enumerate}
		\item $t_0(x_k)\le t_k+t<t_f(x_k)$ for all $t\in [-k,k]$ and all $k$;
		\item there exist a compact set $K\subset \X$ and $\varepsilon>0$ such that $x_k(t_k+t)\in K$ and $|x_k(t_k+t)|\ge \varepsilon$ for all 
		$t\in [-k,k]$ and all $k$; 
		\item for almost all $t\in \R$
		\begin{align} \label{eq:null-o}
		\lim_{k\to \infty} h(t_k+t,x_k(t_k+t),\sigma_k(t_k+t))=0. 
		\end{align}
	\end{enumerate}
\end{defin}
The following definition is a straightforward extension of that given in \cite[Defn.2]{lee_tan_mareels_tac17} for $\X=\R^n$.
\begin{defin} \label{def:wzsd} \rm  Let $f$ and $h$ be as in (\ref{eq:ss}) and (\ref{eq:so}) respectively and let $\T$ be a family of trajectories of (\ref{eq:ss}). The pair $(h,f)$ is weak zero-state detectable (WZSD) with respect to (w.r.t.)  $\T$ if there are no $\T$-zeroing-output sequences for (\ref{eq:ss})-(\ref{eq:so}). 
\end{defin}

The  assumption below plays an important role in checking the UAS of (\ref{eq:ss}) with the help of the output (\ref{eq:so}). It is an extension of condition (H1) in \cite{leejia_tac05} to switched systems. 

\begin{as} \label{ass:integral} \rm There exists a continuous and positive definite function $\alpha:\R_{\ge 0}\to \R_{\ge 0}$ so that for any compact set $K\subset \X$ and any $\mu>0$, there exists $M=M(K,\mu) > 0$ such that for every $(x,\sigma)$ in $\T$ and every $t_0(x_0)\le s<t<t_f(x)$ with $x(\tau)\in K$ for all $\tau\in [s,t]$, we have
\begin{align}
\int_s^t \alpha(|h(\tau,x(\tau),\sigma(\tau))|)\:d\tau \le M+\mu (t-s).
\end{align}
\end{as}

The following result is a reformulation of \cite[Thm.1]{leejia_tac08} that takes into account the fact that the WZSD of a pair $(\tilde h,f)$ w.r.t. a family of trajectories $\T$ implies the output persistent excitation (OPE) of $(\tilde h,f)$ w.r.t. $\T$ (see \cite[Lem.1]{lee_tan_mareels_tac17}; for the definition of the OPE property see \cite{leejia_tac08}). 
\begin{teo} \label{thm:1} \rm Let $\T$ be a US family of trajectories of (\ref{eq:ss}). Suppose that Assumption \ref{ass:integral} holds and that the pair $(h,f)$ is  WZSD w.r.t. $\T$. Then $\T$ is UAS. If, in addition, $\X=\R^n$ and $\T$ is GUS, then $\T$ is GUAS.
\end{teo} 
\begin{IEEEproof} Let $\tilde h=\sqrt{\alpha \circ |h|}$, with $\alpha$ as in Assumption \ref{ass:integral}. Since $(h,f)$ is WZSD w.r.t. $\T$, the same holds for $(\tilde h, f)$. In consequence $(\tilde h,f)$ is OPE w.r.t. $\T$. Hence $\T$ satisfies the hypotheses of Theorem 1 in \cite{leejia_tac08} with $(\tilde h,f)$ in place of $(h,f)$, and then $\T$ is UAS, and GUAS if $\X=\R^n$ and $\T$ is GUS. 
 \end{IEEEproof}

In the applications of Theorem \ref{thm:1} the more difficult task is to check the  WZSD (or the OPE) of the pair $(h,f)$ from their definitions. Although some results for checking that properties were given in \cite{leejia_tac08, lee_tan_nesic_tac15, lee_tan_mareels_tac17} and \cite{lee_tan_mareels_cdc17}, there is still much room for improvement. 

Often the set of trajectories $\T$ satisfies some type of time-dependent constraint, for example its switching signals verify some dwell-time condition, and/or some state-dependent constraint, such as the invariance of $\T$ w.r.t. a closed covering $\chi=\{\chi_i\}_{i\in \I}$ of $\R^n$. We recall that $\T$ is invariant w.r.t. $\chi$ if for all $(x,\sigma)\in \T$, $x(t)\in \chi_{\sigma(t)}$ for all $t\in [t_0(x),t_f(x))$, or, equivalently, $\sigma(t)\in \I_{x(t)}$ for all $t\in [t_0(x),t_f(x))$, where
\begin{align} \label{eq:ixi} 
\I_{\xi}=\{i\in \I:\xi \in \chi_i\}\quad \forall \xi \in \X.
\end{align}
This additional information is, in many cases, useful for checking the WZSD of $(h,f)$ w.r.t. $\T$. 

The precedent discussion motivates the search of solutions for the problems that we face in the paper, and that we briefly state as follows. Given a family $\T$ of trajectories of (\ref{eq:ss}):
\begin{itemize}
\item Obtain easier to check sufficient conditions for the WZSD of a pair $(h,f)$ w.r.t. $\T$, which take explicitly into account the time/state-dependent constraints that $\T$ satisfies (see Section \ref{sec:wzsd}).
\item Give both sufficient and necessary conditions for the UAS or GUAS of $\T$ based on the sufficient conditions for the WZSD obtained (see Sections \ref{sec:guas} and \ref{sec:criteria}).
\end{itemize}
\section{WZSD of switched systems via reduced limiting control systems}\label{sec:wzsd}
In this section we study the WZSD of the pair $(h,f)$ w.r.t. a given family of trajectories $\T$ of (\ref{eq:ss}) through the asymptotic behavior of the $\T$-zeroing-output sequences. In order to motivate the theory we will subsequently develop, we start the section with a motivating example.
\subsection{Motivating example}
Consider the time-invariant switched system (\ref{eq:ss}) with $\X=\R^2$, $\I=\{1,2\}$ and 
\begin{align*}
f_1(\xi)=\left [ \begin{array}{cc} \xi_2 \\
                 -\xi_1
                 \end{array} \right ] 
                 \quad \text{and} \quad
     f_2(\xi)=\left [ \begin{array}{cc} -\xi_1^{1/3}+ a\xi_2 \\
                 -a\xi_1
                 \end{array} \right ],
\end{align*}
with $a>0$. Since $\dot{x}=f_1(x)$ is GUS but not GUAS, the switched system is not GUAS for arbitrary switching. Hence, for obtaining GUAS some restriction on the switched signals must be imposed. For given $T_0>0$ and $\delta_0>0$ consider the set $\S_{T_0,\delta_0}$ of all the switching signals $\sigma$ which satisfy the condition
\begin{align}
 |\{s\in [t,t+T_0]:\sigma(s)=2\}|\ge \delta_0\quad \forall t\ge 0,
\end{align}
that is, the total amount of time the subsystem $2$ is active on $[t,t+T_0]$ is at least $\delta_0$ for all $t\ge 0$.  
We claim that the set $\T$ of trajectories $(x,\sigma)$ with $\sigma \in \S_{T_0,\delta_0}$ is GUAS.

$\T$ is GUS, since $V(\xi)=|\xi|^2/2$ is a common Lyapunov function for (\ref{eq:ss}) because $\partial V/\partial \xi \: f_1(\xi)=0$ and $\partial V/\partial \xi\: f_2(\xi)=-\xi_1^{4/3}\le 0$. Consider the output (\ref{eq:so}) with the time-independent map $h$ defined by $h_1(\xi):= 0$ and $ h_2(\xi):= |\xi_1|$. Since for every $(x,\sigma)\in \T$ and every $t_0(x)\le s<t\le t_f(x)$ 
$$ \int_{s}^t|h(x(\tau),\sigma(\tau))|^{4/3}\:d\tau\le V(x(s))-V(x(t))\le V(x(s)),$$ 
$\T$ satisfies Assumption \ref{ass:integral} with $\alpha(s):= s^{4/3}$ and $M(K,\mu)=\max_{\xi\in K}V(\xi)$.

For fulfilling the hypotheses of Theorem \ref{thm:1} it only remains to check that $(h,f)$ is WZSD w.r.t. $\T$. For a contradiction, suppose there exists a $\T$-zeroing-output sequence \footnote{Recall Definition \ref{def:zeroing-sequence}} $\{(t_k,x_k,\sigma_k)\}$.  Let $z_k(\cdot)=x_k(\cdot+t_k)$ and $\hat{\sigma}_k(\cdot)=\sigma_k(\cdot+t_k)$. From 1) and 2) of Definition \ref{def:zeroing-sequence} it follows that for all $T>0$, the functions $z_k$ are defined on $[-T,T]$ for $k$ large enough and that the sequence $\{z_k\}$ is uniformly bounded on that interval. The latter, the fact that $z_k$ satisfies the differential equation 
\begin{align} \label{eq:zk} 
\dot{z}_k=f(z_k,\hat{\sigma}_k)
 \end{align}
 and Arzela-Ascoli theorem ensure the existence of a subsequence of $\{z_k\}$, which we still denote by $\{z_k\}$, that converges to some continuous function $\bar x:\R\to \R^2$ uniformly on $[-T,T]$ for all $T>0$ ($\bar x$ is a limiting solution of $\T$ according to \cite[Defn.2]{leejia_tac08}). Also due to 2) of Definition \ref{def:zeroing-sequence}, there exists $\varepsilon>0$ such that $|\bar x(t)|\ge \varepsilon$ for all $t\in \R$. From (\ref{eq:zk}) and the uniform convergence of $\{z_k\}$ to $\bar x$ on compact intervals of $\R$, we have that $\bar x$ satisfies the integral equation
\begin{align} \label{eq:lim_int}
 \bar x(t)=\bar x(0)+\lim_{k\to \infty} \int_0^t f(\bar x(s),\hat{\sigma}_k(s))\:ds\quad \forall t\in \R,
\end{align}
while from 3) of Definition \ref{def:zeroing-sequence} and the aforementioned uniform convergence we have that 
\begin{align} \label{eq:lim_output}
\lim_{k\to \infty} h(\bar x(t),\hat{\sigma}_k(t))=0\quad \text{a.e.\;on}\;\R. 
\end{align}
The definition of $h$ and (\ref{eq:lim_output}) imply that 
\begin{align} \label{eq:0-z1}
\bar x_1(t)=0\quad \text{a.e.\;on}\;\mathcal{N},
\end{align}
where $\mathcal{N}=\{t\in \R:\hat \sigma_k(t)=2\;\makebox{for infinitely many}\;k\}$.

When $a=1$, from (\ref{eq:lim_int}), (\ref{eq:0-z1}) and the fact that $f_1(\bar x(t))=f_2(\bar x(t))$ for all $t\in \mathcal{N}$, it follows that $\bar x$ satisfies the differential equation $\dot{\bar x}=f_1(\bar x)$. The latter together with the fact that $\mathcal{N}\cap [0,T_0]$ has infinitely many points imply $\bar x(t)=0$ for all $t\in \R$, which yields a contradiction. This case corresponds to the situation treated in \cite{lee_tan_mareels_tac17}, namely the switched system admits a CZOS, which is $\dot{\bar x}=f_1(\bar x)$.

When $a\neq 1$, there is no CZOS and to draw some conclusion about the behavior of $\bar x$ seems very difficult. For treating this case, we associate to each switching signal $\hat{\sigma}_k$ the function $v_k=[v_{k,1}\; v_{k,2}]':\R\to \R^2$ defined by $v_k(t)=[1\;0]'=e_1$ if $\hat{\sigma}_k(t)=1$ and $v_k(t)=[0\;1]'=e_2$ if $\hat{\sigma}_k(t)=2$. We have that (\ref{eq:lim_int}) is equivalent to
\begin{align} \label{eq:lim_intu}
 \bar x(t)=\bar x(0)+\lim_{k\to \infty} \int_0^t \sum_{i=1}^2 v_{k,i}(s)f_i(\bar x(s))\:ds
\end{align}
and (\ref{eq:lim_output}) is equivalent to 
\begin{align} \label{eq:lim_outputbis}
 \lim_{k\to \infty} \sum_{i=1}^2v_{k,i}(t)h_i(\bar x(t))=0\quad \text{a.e.\;on}\;\R.
\end{align}
Since $v_k$ is piecewise constant and for all $t\in \R$ and all $k$ $v_k(t)\in U$, where $U$ is the convex set
$$U=\{\mu\in \R_{\ge 0}^2:\mu_1+\mu_2=1\}={\rm co}(\{e_1,e_2\}),$$ well-known results of real analysis ensure the existence of a subsequence of $\{v_k\}$, which we still denote by $\{v_k\}$, which weakly converges to some Lebesgue measurable function $\bar u:\R\to U$ (see Definition \ref{def:wc} and Lemma \ref{lem:metric} in Section \ref{sub:limiting}). From the weak convergence of $\{v_k\}$ to $\bar u$ and (\ref{eq:lim_intu}) it follows that for all $t\in \R$,
 $\bar x(t)=\bar x(0)+\int_0^t \sum_{i=1}^2 \bar u_{i}(s)f_i(\bar x(s))\:ds$,
or, equivalently, that $\bar x$ satisfies the differential equation
\begin{align}\label{eq:diff}
 \dot{\bar x}=\bar u_1f_1(\bar x)+\bar u_2 f_2(\bar x)
\end{align}
on $\R$. From the weak convergence of $\{v_k\}$ to $\bar u$, (\ref{eq:lim_outputbis}) and the Lebesgue convergence theorem, we have that for all $T>0$
\begin{align}
 \int_{-T}^{T} \sum_{i=1}^2 \bar u_{i}(t)h_i(\bar x(t))\:dt=0. \label{eq:0-int}
\end{align}
The nonnegativeness of the integrand in (\ref{eq:0-int}) and the arbitrariness of $T$ imply that
\begin{align} \label{eq:0-out}
 \bar y(t)=\bar u_1(t)h_1(\bar x(t))+\bar u_2(t)h_2(\bar x(t))=0 \;\text{a.e.\;on}\;\R.
\end{align}
From (\ref{eq:diff}) and (\ref{eq:0-out}) we have that the pair $(\bar x,\bar u)$ satisfies 
\begin{align} \label{eq:rcsex} 
\left \{ \begin{array}{rcl} \dot{\bar x}_1(t) & = &[\bar u_1(t)+a\bar u_2(t)]\bar x_2(t)\\
 \dot{\bar x}_2(t) & = & -\bar u_1(t)\bar x_1(t)\\
 \bar y(t) & = & \bar u_2(t) |\bar x_1(t)|=0 \end{array} \right .,
\end{align}
for almost all $t\in \R$. Here we have dropped the terms containing $\bar x_1(t)$ in $f_2(\bar x(t))$ since they vanish when $u_2(t)\neq 0$.

Since the switching signals $\sigma_k$ belong to $\S_{T_0,\delta_0}$, $\delta_0\le |\{s\in [0, T_0]:\hat \sigma_k(s)=2\}|=\int_0^{T_0} v_{k,2}(t)\:dt$, which together with the weak convergence of $\{v_k\}$ to $\{\bar u\}$ imply that
 $\int_0^{T_0}\bar u_2(t)\:dt\ge \delta_0$
and then that
\begin{align}
\label{eq:measure}
 |\{t\in [0,T_0]:\bar u_2(t)>0\}|\ge \delta_0.
\end{align}

Let $F\subset [0,T_0]$ be the set where the equalities in (\ref{eq:rcsex}) hold and let $E=\{t\in F:\bar u_2(t)>0\}$. Then $|E|\ge \delta_0$. Note that $\bar x_1(t)=0$ for all $t\in E$. To continue with our analysis we need the following result, whose proof is given in the Appendix.
\begin{lema} \label{lem:0-der}
 Let $\varphi:(a,b)\to \R$ and let $E\subset (a,b)$ be a subset of positive measure. Suppose that for all $t\in E$, $\varphi$ is differentiable at $t$ and $\varphi(t)=0$. Then there exists a subset $E^*\subset E$ such that $|E \setminus E^*|=0$ and $\dot{\varphi}(t)=0$ for all $t \in E^*$.  
 \end{lema}
By applying Lemma \ref{lem:0-der} to $\bar x_1$, we have that $\dot{\bar x}_1(t)=0$ for all $t\in E^*$, where $E^*$ is a certain nonempty subset of $E$. Then, from the first equation in (\ref{eq:rcsex}) and that fact that $\bar u_1(t)+2\bar u_2(t)>0$ it follows that $\bar x_2(t)=0$ for all $t\in E^*$. Consequently $\bar x(t)=0$ on $E^*$, which contradicts the fact that $|\bar x(t)|\ge \varepsilon$ for all $t\in \R$. Therefore, the pair $(h,f)$ is WZSD w.r.t. $\T$ and the family of trajectories is GUAS. 

We note that the stability analysis of $\T$ is challenging since one of the subsystems is only GUS and the set of switching signals $\S_{T_0,\delta_0}$ contains switching signals $\sigma$ for which the corresponding sequence of switching times $\{\tau_k^{\sigma}\}$ satisfies the condition $\tau^{\sigma}_{k+1}-\tau^{\sigma}_k\to 0$. Up to our best knowledge no existing GUAS criteria can be straightforwardly applied to this case. 
On the other hand, it also seems diffilcult in this case to apply the tools given in \cite{lee_tan_nesic_tac15} and in \cite{lee_tan_mareels_cdc17} for checking the OPE condition.     
\begin{remark}
Note that the key facts which allow us to prove the WZSD of the pair $(h,f)$ are: i) for every $\T$-zeroing-output sequence $\{(t_k,x_k,\sigma_k)\}$ the limits $(\bar x,\bar u)$ of ``convergent" subsequences of $\{(z_k,\hat \sigma_k)\}$ are solutions of the ``reduced control system" (\ref{eq:rcsex}) with controls $\bar u$ restricted to satisfy $\bar u(t)\in U$ for all $t\in \R$ and the constraint (\ref{eq:measure}) induced by the restriction imposed to the switching signals; and ii) for every pair $(\bar x,\bar u)$ satisfying (\ref{eq:rcsex}), $\bar u(t)\in U$ for all $t\in \R$  and (\ref{eq:measure}), $\inf_{t \in \R}|\bar x(t)|=0$, i.e. the reduced control system (\ref{eq:rcsex}) enjoys some kind of detectability property. \mer
\end{remark}
In the following subsections we formalize and extend the procedure used in the motivating  example to switched NLTV systems obtaining in that way sufficient conditions for the WZSD of a pair $(h,f)$ w.r.t. a family of trajectories. These conditions  are in many cases easier to check than the existing ones in present literature.  
                 
\subsection{Embedding of the switched system into a control one}
For characterizing the asymptotic behavior of the $\T$-zeroing-output sequences of (\ref{eq:ss})-(\ref{eq:so}) and in this way to study the WZSD of a pair $(h,f)$ w.r.t. the a family $\T$ as in the motivating example, it is convenient to embed system (\ref{eq:ss})-(\ref{eq:so}) into the control-affine time-varying system with outputs
\begin{eqnarray}
 \dot{x}  =\sum_{i=1}^Nu_i f_i(t,x):=F(t,x) u \label{eq:cs}\\
y =\sum_{i=1}^Nu_ih_i(t,x):=H(t,x)u \label{eq:coutput}
\end{eqnarray}
where, $u=[u_1\cdots u_N]'$ and, for all $t\ge 0$ and all $\xi\in \X$, 
$F(t,\xi):=[f_1(t,\xi)\ldots f_N(t,\xi)]\in \R^{n\times
N}$ and $H(t,\xi):=[h_1(t,\xi)\ldots h_N(t,\xi)]\in \R^{p\times
N}$. We assume that the admissible controls $u$ of (\ref{eq:cs}) belong to $\U$, the set of all the Lebesgue measurable functions $u:\R \to U$, where $U$ is the convex set
$$U=\left \{\mu \in \R_{\ge 0}^N:\sum_{i=1}^N\mu_i=1 \right \}.$$
Note that $U={\rm co}(U^*)$, with $U^*:=\{e_1,\ldots,e_N\}$, and where $e_i \in \R^N$ denotes the $i$-th canonical vector of $\R^N$. 

The embedding of (\ref{eq:ss}) into (\ref{eq:cs}) is performed by identifying the set $\S$ of all the switching signals with the set $\U_{pc}^*$ of all controls $u\in \U$ that take values in $U^*$ and are piecewise constant and continuous from the right, by means of the bijection $\sigma \mapsto u_\sigma$, $u_\sigma(\cdot)=e_{\sigma(\cdot)}$. Since $f(t,\xi,i)\equiv F(t,\xi)e_i$, and $h(t,\xi,i)\equiv H(t,\xi)e_i$, the solutions of (\ref{eq:ss}) corresponding to a switching signal $\sigma$ are the same as those of (\ref{eq:cs}) which correspond to the control $u_\sigma$, and for any solution $x:[t_0(x),t_f(x))\to \R^n$ of (\ref{eq:ss}) corresponding to a switching signal $\sigma$, 
\begin{align} \label{eq:equalout}
h(t, x(t),\sigma(t))=H(t,x(t))u_{\sigma}(t)\quad \forall t \in [t_0(x),t_f(x)).
\end{align}
From now on we will identify every switching signal $\sigma$ with the corresponding control $u_{\sigma}$, and write $\sigma$ in place of $u_{\sigma}$ without risk of confusion. In other words, depending on the context, $\sigma$ will represent either a switching signal or the corresponding piecewise control $u_{\sigma}$. We will also identify the set $\S$ with $\U^*_{pc}$ and any family $\T$ of trajectories $(x,\sigma)$ of (\ref{eq:ss})-(\ref{eq:so}) with the corresponding family of trajectories $(x,u_{\sigma})$ of (\ref{eq:cs})-(\ref{eq:coutput}). 
\subsection{Limiting trajectories} \label{sub:limiting}
The concept of limiting trajectory of $\T$ that we introduce next, which is inspired in that of limiting solution for time-varying ordinary differential equations (see \cite{Artstein76}), is useful for studying the asymptotic behavior of sequences of trajectories of (\ref{eq:ss}). Its definition requires the following standard notion of weak-convergence in $\U$.
\begin{defin}\label{def:wc}
Given a sequence $\{u_k\}$ in $\U$ and $u \in \U$, we say that $u_k\rightharpoonup u$ if for all 
$f\in L^1_N(\R)$
$$ \lim_{k \to \infty} \int_{-\infty}^{\infty} f(t)'u_k(t) \: dt =\int_{-\infty}^\infty f(t)'u(t) \:dt. $$

Note that if $u_k \rightharpoonup u$, then for every $a < b$ 
$$ \lim_{k \to \infty} \int_{a}^{b} f(t)'u_k(t) \: dt =\int_{a}^b f(t)'u(t) \:dt, $$
for all $f\in L^1_N([a,b])$. \mer
 \end{defin}
Next we introduce the definitions of limiting trajectory of a set $\T$ of trajectories of the switched system (\ref{eq:ss}) and of $\T$-zeroing-output limiting trajectory of (\ref{eq:ss})-(\ref{eq:so}). 
\begin{defin} \label{def:limsol} \rm A pair $(\bar x,\bar u)$, with $\bar x:\R \to \X$ and $\bar u \in \U$ is a limiting trajectory of $\T$ corresponding to a sequence $\gamma=\{t_k\}$ in $\R_{\ge 0}$ such that $t_k\to \infty$, if there exist a sequence $\{(x_k,\sigma_k)\}$ in $\T$ and a compact set  $K\subset \X$ such that:
 	\begin{enumerate}
 		\item $x_k(t)\in K$ 
 		for all $t\in [t_k-k,t_k+k]\subset [t_0(x_k),t_f(x_k))$ and all $k$, and  $\{x_k(t_k+\cdot)\}$
 		converges to $\bar x$ uniformly on $[-T,T]$ for all $T>0$; and 
 		\item $\sigma_k(t_k+\cdot)\rightharpoonup \bar u$.
 	\end{enumerate}
 	If, in addition, $\{(t_k,x_k,\sigma_k)\}$ is a $\T$-zeroing-output sequence of (\ref{eq:ss})-(\ref{eq:so}), then we say that $(\bar x,\bar u)$ is a $\T$-zeroing-output limiting trajectory of (\ref{eq:ss})-(\ref{eq:so}).
 \end{defin}
\begin{remark} 
For each limiting trajectory $(\bar x,\bar u)$ of $\T$, $\bar x$ is a limiting solution of $\T$ in the sense of \cite[Defn.2]{leejia_tac08}. The notion of limiting solution for switched systems was introduced in \cite{leejia_tac08} in order to study the asymptotic behavior of solutions of the switched system and used in \cite{lee_tan_mareels_tac17} for studying the WZSD of a pair $(h,f)$. As a difference with those works, the consideration of limiting trajectories instead of just limiting solutions, will allow us to characterize the $\T$-zeroing-output limiting trajectories as solutions of certain control systems. This characterization will ease the analysis of the WZSD of a pair $(h,f)$. \mer
\end{remark} 

For guaranteeing  the existence of limiting trajectories we assume that $f$ satisfies the following boundedness condition.  
\begin{as}\label{ass:boundedness} \rm $f_i$ is uniformly bounded for every $i\in \I$.
\end{as}
\begin{lema} \rm \label{lm:subseq} Let $f$ in (\ref{eq:ss}) satisfy Assumption \ref{ass:boundedness}. Let $\gamma=\{t_k\}$ be a sequence of positive real numbers such that $t_k\to \infty$ and let $\{(x_k,\sigma_k)\}$ be a sequence of trajectories of (\ref{eq:ss}) such that for some compact 
	set $K\subset \X$ and for every $k$, $x_k(t)\in  K$ for all $t\in [t_k-k,t_k+k]\subset [t_0(x_k),t_f(x_k))$. Then there exist a 
	subsequence $\{(x_{k_l},\sigma_{k_l})\}$, a continuous function $\bar x:\R\to \X$ and a control $\bar u\in \U$ such that 
	$\{x_{k_l}(t_{k_l}+\cdot)\}$ converges to $\bar x$ uniformly on $[-T,T]$ for all $T>0$ and $\sigma_{k_l}(t_{k_l}+\cdot)\rightharpoonup \bar u$.
\end{lema}

Lemma \ref{lm:subseq} easily follows from the Arzela-Ascoli theorem and the following result about the sequential compactness of $\U$.
\begin{lema} \rm
   \label{lem:metric}
   For every sequence $\{u_k\}$ in $\U$ there exist $u\in \U$ and a subsequence $\{u_{k_l}\}$ such that $u_{k_l}\rightharpoonup  u$.
\end{lema}

Lemma \ref{lem:metric} is a straightforward consequence of the following  facts: i) the unitary ball $\mathcal{B}_1^{\infty}$ of $L^{\infty}_{N}(\R)$ endowed with the weak$^*$ topology (as the dual space of $L^1_N(\R)$) is a compact metric space, due to Alaoglu's theorem and the separability of $L^1_N(\R)$ \cite[Thm.3.16]{rudin_functional_analysis}, and ii) $\U$ is a weak$^*$ closed subset of $\mathcal{B}_1^{\infty}$. 

When $\T$ is invariant w.r.t. a closed covering $\chi$,  the limiting trajectories of $\T$ satisfy an analogous state-dependent constraint. Define for $\xi \in \R^n$ the following set of control values:
\begin{align}
U_{\xi}={\rm co}\{e_i:i\in \I_{\xi}\},
\end{align}
where $\I_{\xi}$ is given by (\ref{eq:ixi}).

The next lemma, whose proof is given in the Appendix, establishes the aforementioned state-dependent constraint.

\begin{lema}\label{lem:sdc} \rm Let $\T$ be a family of trajectories of $(\ref{eq:ss})$ which is invariant w.r.t. a closed covering $\chi$. Let $(\bar x, \bar u)$ be a limiting trajectory of $\T$. Then $\bar u(t)\in U_{\bar x(t)}$ for almost all $t\in \R$.  
 \end{lema}
\subsection{Reduced limiting control systems}

In this subsection we introduce the notion of reduced limiting control systems for a family of trajectories $\T$ of (\ref{eq:ss}). Roughly speaking, these systems are auxiliary control systems with outputs, such that each $\T$-zeroing-output limiting trajectory $(\bar x,\bar u)$ is a trajectory of some of them for which the output is identically zero and the control $\bar u$ satisfies certain constraints induced by those that the switching signals satisfy. Hence, the behavior of the trajectories of these auxiliary control systems will give us information about the $\T$-zeroing-output limiting trajectories, and consequently about the WZSD of the pair $(h,f)$.

 In the following we assume that $\T$ is invariant w.r.t. some closed covering $\chi$. This assumption does not imply any loss of generality since $\T$ is always invariant w.r.t. the trivial covering $\chi=\{\chi_i\}_{i=1}^N$, where $\chi_i=\R^n$ for all $i$. 

For defining the reduced limiting control systems associated to the family $\T$ we need to make further assumptions on $f$ and $h$. One of these assumptions is concerned with the asymptotic behavior of the time-translations of certain functions related to $f$ and $h$. With this aim we recall some definitions. 

A sequence of non negative times $\gamma=\{t_k\}$ is said to be admissible for a function  $g:\R_{\ge 0}\times \X \to \R^m$, if $t_k\to \infty$ and there exist a function $g_{\gamma}:\R\times \X \to \R^m$ and a  set $E_{\gamma}\subset \R$ of Lebesgue measure zero such that
$$ \lim_{k \to \infty}g(t_k+s,\xi)=g_{\gamma}(s,\xi)\quad \forall \xi \in \R^n,\;\forall s\in \R\setminus E_{\gamma}.$$
The function $g_\gamma$ is called the limiting function of $g$ associated to $\gamma$. We will denote by $L(g)$ the set of all the limiting functions of $g$ and by $\Lambda(g)$ the set of its admissible sequences. We say that $g$ is precompact if for any sequence $\{t_k\}$ in $\R_{\ge 0}$ such that $t_k\to \infty$ there exists a subsequence $\gamma$ of $\{t_k\}$ such that $\gamma\in \Lambda(g)$. If $g(t,\xi)$ is measurable in $t$ and continuous in $\xi$, then every $g_{\gamma}\in L(g)$ is measurable in $t$, but it is not necessarily continuous in $\xi$. For guaranteeing the continuity of $g_\gamma$ in $\xi$, $g$ has to verify stronger continuity conditions. It can be easily proved that if $g(t,\xi)$ is continuous in $\xi$, uniformly in $t$, then the same holds for every $g_{\gamma}\in L(g)$.
\begin{remark} \rm The notion of precompact function we consider is stronger than that considered in \cite{Artstein76} for defining the limiting equations of a time-varying system of ordinary differential equations. Such a stronger definition is needed in our framework for dealing with the type of convergence we consider in the space of controls, since the notion used in \cite{Artstein76} is too weak in this case.  \mer
\end{remark}
\begin{remark} \rm \label{rem:aap-precompact}
 A class of functions which is precompact in our sense is that of almost asymptotic periodic (AAP) ones (see \cite{leejia_tac05} for its definition and some interesting properties). That class includes that of continuous functions which do not depend on $t$. 
 \mer
\end{remark}
We will also employ the concept of zeroing pair (\cite{lee_tan_nesic_tac15}) in order to relax regularity requirements on the functions $f$ and $h$ and simplify the form of the reduced limiting control systems. Let $g:\R_{\ge 0}\times \X\to \R^m$ and $\hat g:\R_{\ge 0}\times \X\to \R^{\hat q}$. The pair $(g,\hat g)$ is a zeroing pair if for any sequence $\{t_k\}$ in $\R_{\ge 0}$ with $t_k\to \infty$, any compact subset $K\subset \X$, any constant $\varepsilon>0$ and any sequence $\{v_k\}$ in $K$ such that $|v_k|\ge \varepsilon$ for all $k$, the following implication holds: $g(t_k,v_k)\to 0 \; \Rightarrow \; \hat g(t_k,v_k)\to 0.$

Now we are in position to state the assumptions we make on $f$ and $h$ in order to define the reduced limiting control systems associated to $\T$.
\begin{as} \label{ass:aap}\rm $f$ and $h$ in (\ref{eq:ss})-(\ref{eq:so}) satisfy for all $i\in \I$ the following:
\begin{enumerate}	
\item \label{item:h} $h_i(t,\xi)$ is uniformly bounded, measurable in $t$, continuous in $\xi$, uniformly in $t$, and precompact;
\item  \label{eq:decomp} $f_i$ can be decomposed as $f_i=\hat{f}_i+\Delta f_i$, with
\begin{enumerate}
	\item $\hat{f}_i(t,\xi)$  uniformly bounded, measurable in $t$, continuous in $\xi$, uniformly in $t$, and precompact;
    \item $\Delta f_i$  such that $(h_i,\Delta f_i)$ is a zeroing pairs.
\end{enumerate}
\end{enumerate}
\end{as}
We next provide some comments about the required Assumption \ref{ass:aap}. First we note that, loosely speaking, condition \ref{eq:decomp}) requires that $f_i$ can be decomposed as the sum of two terms, one of them, $\Delta f_i$, converging to zero as $h_i\to 0$ and $t\to \infty$. Such a decomposition is always possible, since we can take $\Delta f_i\equiv 0$. Nevertheless, the possibility of selecting $\Delta f_i$ in other ways will allow us to relax the regularity requirements on $f_i$ and to simplify the form of the reduced limiting systems. With respect to the regularity conditions 1) and 2.a) imposed on $h_i$ and on $\hat{f}_i$  respectively, they are satisfied by the class of AAP functions, since any AAP function is measurable in $t$, continuous in $\xi$, uniformly in $t$, and precompact (see \cite[Sec.III]{leejia_tac05}). Hence, 1) and 2.a) hold when $h_i$ and $\hat f_i$ do not depend on $t$. In particular, Assumption \ref{ass:aap} is indeed guaranteed to hold when the switched system (\ref{eq:ss})-(\ref{eq:so}) is time-invariant.

For a system (\ref{eq:ss}) that satisfies Assumption \ref{ass:aap}, we define the notion of reduced limiting control system for a set of trajectories $\T$ of (\ref{eq:ss}) which is invariant w.r.t. a closed covering as follows. Let $\S_{\T}=\{\sigma \in \S:\exists\, x \;\text{such that}\;(x,\sigma)\in \T\}$ and let $\S^*_{\T}$ be the set of controls $u\in \U$ for which there exist $\{t_k\}$ such that $t_k\to \infty$ and $\{\sigma_k\} \in \S_{\T}$ such that $\sigma_{k}(t_k+\cdot)\rightharpoonup  u$.

For any $\gamma \in \Lambda:=\cap_{i\in \I}[\Lambda(\hat{f}_i) \cap \Lambda(h_i)]$, let, for all $t\in \R$ and all $\xi \in \R^n$, 
$\hat{F}_{\gamma}(t,\xi):=[\hat{f}_{1,\gamma}(t,\xi)\ldots \hat{f}_{N,\gamma}(t,\xi)]$ and 
$\hat H_{\gamma}(t,\xi):=[|h_{1,\gamma}(t,\xi)|\ldots |h_{N,\gamma}(t,\xi)|]$. Then 
\begin{align} 
\label{eq:lcs}
\Sigma_{\gamma}: \;\left \{ \begin{array}{l} \dot{x} = \hat{F }_{\gamma}(t, x) u  \\
y = \hat H_{\gamma}(t,x)u=0
\end{array}
\right .,\; u(t)\in U_{x(t)},\;u \in \S^*_{\T}.
\end{align}
is a reduced limiting control system for $\T$. 

We say that $(x,u)$ is a complete trajectory of $\Sigma_{\gamma}$ if $x:\R\to \X$ is locally absolutely continuous, $u\in \S^*_{\T}$ and, for almost all $t\in \R$, $u(t)\in U_{x(t)}$, $\dot{x}(t)=\hat{F}_{\gamma}(t,x(t))u(t)$ and $y(t)=\hat H_{\gamma}(t,x(t))u(t)=0$. We also say that such a complete trajectory is bounded if there exists a compact subset $K\subset \X$ such that $x(t)\in K$ for all $t\in \R$.
\begin{remark}\rm The definition of reduced limiting control system is inspired --and extends to switched systems-- that of reduced limiting system for
 ordinary differential equations introduced in \cite{leejia_tac05}. In fact, when (\ref{eq:ss})-(\ref{eq:so}) is a nonswitched system, {\em i.e.}, $\dot{x}=f(t,x)$ and $y=h(t,x)$, the systems $\Sigma_{\gamma}$ coincide with the reduced limiting systems defined in \cite{leejia_tac05}.   \mer
 \end{remark}
We next provide some discussion on the reduced limiting control systems just introduced. When the functions $\hat{f}_i$ and $h_i$ in Assumption \ref{ass:aap} are time-independent, that is $\hat f_i(t,\xi)\equiv \bar f_i(\xi)$ and $h_i(t,\xi)\equiv \bar h_i(\xi)$, there is only one reduced limiting control system, which has the form 
\begin{align*} 
\Sigma: \;\left \{ \begin{array}{l} \dot{x} = \sum_{i=i}^N u_i \bar f_i(x)\\
y = \sum_{i=i}^N  u_i |\bar h_i(x)| =0
\end{array}
\right .,\; u(t)\in U_{x(t)},\;u \in \S^*_{\T}.
\end{align*}
This is the case of time-invariant switched systems or, more generally, time-invariant switched systems perturbed by additive time-varying terms depending on the state through the output that converges to zero as $t\to \infty$ and the output converges to zero, i.e. $f_i(t,\xi)=\bar f_i(\xi)+\omega_i(t,\bar h_i(\xi))$, with $\omega(t,y)\to 0$ as $t\to \infty$ and $y\to 0$. When the functions $\hat f_i$ and $\hat h_i$ are time-dependent, it may be very difficult or even impossible to find the limiting functions $\hat{f}_{i,\gamma}$ and $h_{i,\gamma}$. Nevertheless, in some cases it is not necessary to know the exact form of the limiting functions for performing the analysis of its solutions, as we will see in Examples \ref{ex:1} and \ref{ex:3}, where the only thing that we need to know is the nonzero property of the limiting functions.  

The reduced limiting control systems $\Sigma_{\gamma}$ incorporate the time/state-dependent constraints that the family of trajectories $\T$ satisfies through the restrictions the controls have to satisfy. The invariance of $\T$ w.r.t. the covering $\chi$ is taken into account through the restriction $u(t)\in U_{x(t)}$. The time-dependent constraints the switching signals of the trajectories of $\T$ satisfy are reflected in the restriction $\sigma \in \S^*_{\T}$. For example, if the switching signals in $\S_{\T}$ have common average dwell-time $\tau_D>0$ and chattering bound $N_0\in \N$ (see \cite{Liberzonbook}), then the controls in $\S^*_{\T}$ are controls in $\U^*_{pc}$ (and therefore switching signals) which have average dwell-time $\tau_D>0$ and chattering bound $N_0\in \N$. More generally, it can be proved that if $\S_{\T}$ is contained in a set $\mathcal{V}\subset \U_{pc}^*$ which is invariant for time-translations and sequentially compact with respect to the almost everywhere convergence then $\S^*_{\T}\subset \mathcal{V}$ (see \cite{mancilla-garcia-nolcos, mancilla-garcia-2} for examples of such  sets $\mathcal{V}$). In many practical cases we do not need to know exactly $\S_{\T}^*$, but some useful property the controls in $\S_{\T}^*$ enjoy. This is case in the Motivating Example, where knowing that the controls satisfy the constraint (\ref{eq:measure}) was enough for drawing a conclusion (see also Example \ref{ex:3} in Section \ref{ej}).   

The following result establishes the connection between the $\T$-zeroing-output limiting trajectories of (\ref{eq:ss})-(\ref{eq:so}) and the family $\{\Sigma_{\gamma}\}_{\gamma \in \Lambda}$ of reduced limiting control systems for $\T$.
\begin{teo} \label{teo:limiting} \rm Suppose that $f$ and $h$ in (\ref{eq:ss})-(\ref{eq:so}) satisfy Assumptions \ref{ass:boundedness} and \ref{ass:aap} and let $\T$ be a family of trajectories of (\ref{eq:ss}) which is invariant w.r.t. a closed covering $\chi$. If $(\bar x,\bar u)$ is a $\T$-zeroing-output limiting trajectory of (\ref{eq:ss})-(\ref{eq:so}), then there exists $\gamma \in \Lambda$ such that $(\bar x,\bar u)$ is a bounded complete trajectory of the reduced limiting control system $\Sigma_\gamma$. 
 \end{teo}
 \begin{remark} Theorem \ref{teo:limiting} states that the possible asymptotic behaviors of $\T$-zeroing-output sequences of (\ref{eq:ss})-(\ref{eq:so}), which are represented by the $\T$-zeroing-output limiting trajectories, can be studied by analyzing, for each $\gamma \in \Lambda$, the behavior of the bounded complete trajectories $(x,u)$ of the system $\Sigma_\gamma$, i.e., by analyzing the behavior of the solutions of reduced systems of differential equations. \mer
 \end{remark}
\begin{IEEEproof}[Proof of Theorem~\ref{teo:limiting}]  Let $(\bar x,\bar u)$ be a $\T$-zeroing-output limiting trajectory of (\ref{eq:ss})-(\ref{eq:so}). Then there exists a $\T$-zeroing-output sequence $\{(t_k,x_k,\sigma_k)\}$ for which 1) and  2) of Definition \ref{def:limsol} hold and, in addition, for almost all $s\in \R$ 
\begin{align}\label{eq:z-output}
  \lim_{k\to \infty}|h_{\sigma(t_k+s)}(t_k+s,x_k(t_k+s))|=0.
\end{align}
Since the functions $\hat{f}_i$ and ${h}_i$ are precompact ones for all $i\in \I$, passing to a subsequence and relabeling if necessary, we can assume that $\gamma=\{t_k\}\in \Lambda$. We will show that $(\bar x, \bar u)$ is a bounded complete trajectory of $\Sigma_{\gamma}$.

From 1) of  Definition \ref{def:limsol} it follows that $\bar x$ is bounded. Since $\sigma_k\in \S_{\T}$ for all $k$ and taking into account 2) of Definition \ref{def:limsol}, it follows that $\bar u \in \S^*_{\T}$. From Lemma \ref{lem:sdc}, we have that $\bar u(t)\in U_{\bar x(t)}$ a.e. on $\R$.

For $k\in \N$, let $F_k(\cdot,\cdot)=F(t_k+\cdot,\cdot)$, $\hat F_k(\cdot,\cdot)=\hat F(t_k+\cdot,\cdot)$ and $\Delta F_k(\cdot,\cdot)=\Delta F(t_k+\cdot,\cdot)$, where the functions $\hat F$ and $\Delta F$ are defined in the same way as $F$ but with $\hat f_i$ and $\Delta f_i$ instead of $f_i$. Since $f_i$ and $\hat f_i$ are uniformly bounded functions due to Assumptions \ref{ass:boundedness} and \ref{ass:aap}, respectively, we have that $F$, $\hat F$ and $\Delta F$ have uniformly bounded components. We also have that the components of $\hat F(t,\xi)$ are continuous in $\xi$, uniformly in $t$, due to Assumption \ref{ass:aap}. 

Let $z_k(\cdot)=x_k(t_k+\cdot)$ and $v_k(\cdot)=\sigma(t_k+\cdot)$ for all $k\in \N$. Note that $z_k\to \bar x$ uniformly on any compact interval and that there exist a compact set $K\subset \X$ and a positive constant $\varepsilon_0$ such that $z_k(t)\in K$ and $|z_k(t)|\ge \varepsilon_0$ for all $t\in [-k,k]$ and all $k$ since $\{(t_k,x_k,\sigma_k)\}$ is a $\T$-zeroing-output sequence for (\ref{eq:ss})-(\ref{eq:so}).  

For any $t \in \R$ and by using Assumption \ref{ass:aap} and the fact that $(x_k,\sigma_k)$ is a trajectory of (\ref{eq:cs}), we have that
\begin{align} \label{eq:lim1}
 \bar x(t) &=\lim_{k\to \infty} \left [ z_k(0)+\int_0^t F_k(s,z_k(s))v_k(s)\:ds\right ] \nonumber\\
 &= \bar x(0) + \lim_{k\to \infty} \left [\int_0^t \hat F_k(s,z_k(s))v_k(s)\:ds\right. \nonumber \\ &\left. \qquad \quad + \int_0^t \Delta F_k(s,z_k(s))v_k(s)\:ds \right ].
\end{align}

Let $s\in \R$ be such that (\ref{eq:z-output}) holds. Then, from the facts: i) $\{z_k(s)\}$ is a sequence in $K$ such that $|z_k(s)|\ge \varepsilon_0$  for all $k$ large enough, ii) $(h_i,\Delta f_i)$ is a zeroing-pair for all $i\in \I$ due to Assumption \ref{ass:aap} and iii)  $\I$ is a finite set, it follows that 
\begin{multline*}
 \lim_{k\to \infty} \Delta F_k(s,z_k(s))v_k(s) = \\ \lim_{k\to \infty} \Delta f_{\sigma(t_k+s)}(t_k+s,z_k(s))=0.
\end{multline*}

In consequence $\lim_{k\to \infty}\Delta F_k(s,z_k(s))v_k(s)=0$ a.e. on $\R$. The uniform boundedness of $\{z_k\}$ and of the components of $\Delta F$ and Lebesgue's convergence theorem imply that
\begin{align} \label{eq:lim2}
 \lim_{k\to \infty} \int_0^t \Delta F_k(s,z_k(s))v_k(s)\:ds=0.
\end{align}

On the other hand, 
\begin{multline*}
\int_0^t \hat F_k(s,z_k(s))v_k(s)\:ds =\\ \int_0^t [\hat F_k(s,z_k(s))-\hat F_k(s,\bar x(s))]v_k(s)\:ds \\
                                                          + \int_0^t [\hat F_{k}(s,\bar x(s))-\hat F_{\gamma}(s,\bar x(s))]v_k(s)\:ds\\
                                                          + \int_0^t \hat F_{\gamma}(s,\bar x(s))(v_k(s)-\bar u(s))\:ds \\
                                                          + \int_0^t \hat F_{\gamma}(s,\bar x(s))\bar u(s)\:ds.
\end{multline*}
Taking into account the uniform convergence of $z_k$ to $\bar x$ on compact intervals, the fact that $z_k$ and $\bar x$ take values in the compact set $K$, the uniform continuity in the second argument of the components of $\hat F$ and the boundedness of the components of $v_k$ we have that
\begin{align*} 
 \lim_{k\to \infty} \int_0^t [\hat F_k(s,z_k(s))-\hat F_{k}(s,\bar x(s))]v_k(s)\:ds =0.
\end{align*}
The convergence of the components of $\hat F_k(s,\bar x(s))$ to that of $\hat F_{\gamma}(s,\bar x(s))$ for almost all $s\in \R$, the fact that $\bar x$ takes values in the compact set $K$, the uniform boundedness of the components of $F$, the boundedness of the components of $v_k$ and the Lebesgue convergence theorem imply that
\begin{align*}
 \lim_{k\to \infty} \int_0^t [\hat F_k(s,\bar x(s))-\hat F_{\gamma}(s,\bar x(s))]v_k(s)\:ds=0.
 \end{align*}
The continuity of $\hat F_{\gamma}(\cdot,\bar x(\cdot))$ and the weak convergence of $v_k$ to $\bar u$ yield
\begin{align*} 
 \lim_{k\to \infty} \int_0^t \hat F_{\gamma}(s,\bar x(s))(v_k(s)-\bar u(s))\:ds=0.
\end{align*}
In consequence
\begin{align} \label{eq:lim3}
 \lim_{k\to \infty} \int_0^t \hat F_k(s,z_k(s))v_k(s)\:ds=\int_0^t \hat F_{\gamma}(s,\bar x(s))\bar u(s)\:ds.
\end{align}
Then, from (\ref{eq:lim1}), (\ref{eq:lim2}) and (\ref{eq:lim3}) it follows that for all $t\in \R$
\begin{align*}
 \bar x(t)= \bar x(0)+\int_0^t \hat F_{\gamma}(s,\bar x(s))\bar u(s)\:ds,
\end{align*}
or, equivalently, that $\dot{\hat x}(t)=\hat F_{\gamma}(s,\bar x(t))\bar u(t)$ a.e. on $\R$. 

In order to show that $\hat H_{\gamma}(t,\bar x(t))\bar u(t)=0$ a.e. on $\R$, we proceed as follows. Let $\hat H=[|h_1|\cdots |h_N|]$ and $\hat H_k(\cdot,\cdot)=\hat H(t_k+\cdot,\cdot)$. Due to Assumption \ref{ass:aap}, the components of $\hat H(t,\xi)$ are uniformly bounded and continuous in $\xi$, uniformly in $t$.  

Let $a$ and $b \in \R$ be such that $a<b$. From the facts that $z_k(s)\in K$ if $k$ is large enough, the uniform boundedness of the components of $\hat H$, the uniform boundedness of the components of $\{v_k\}$, the convergence of $\hat H_k(\cdot,z_k(\cdot))v_k(\cdot)$ to zero a.e. on $\R$ (due to (\ref{eq:z-output})), and Lebesgue's convergence theorem yield
\begin{align*}
  \lim_{k\to \infty}\int_{a}^b \hat H_k(s,z_k(s)) v_k(s)\:ds=0.
\end{align*}

Then, by using arguments similar to those used for proving (\ref{eq:lim3}), one can show that
\begin{align*}
 0=\lim_{k\to \infty}\int_{a}^b \hat H_k(s,z_k(s)) v_k(s)\,ds = \int_a^b \hat H_{\gamma}(s,\bar x(s))\bar u(s)\,ds.
\end{align*}
Finally, the arbitrariness of $a$ and $b$ and the nonnegativeness of $\hat H_{\gamma}(s,\bar x(s))\bar u(s)$  imply that $\hat H_{\gamma}(t,\bar x(t))\bar u(t)=0$ a.e. on $\R$.
\end{IEEEproof}
 \subsection{A criterion for WZSD}
 We next give a criterion to determine the WZSD of a pair $(h,f)$ w.r.t. a family of trajectories $\T$ of (\ref{eq:ss}). This criterion is formulated in terms of the WZSD of the reduced limiting control systems introduced in the precedent subsection.
 \begin{defin}\rm A reduced limiting control system $\Sigma_{\gamma}$ for $\T$ is WZSD if for every bounded complete trajectory $(x,u)$ of $\Sigma_{\gamma}$ 
we have that $\inf_{t\in \R}|x(t)|=0$.
 \end{defin}
\begin{teo}\rm \label{thm:criterio_wzsd} Let Assumptions \ref{ass:boundedness} and \ref{ass:aap} hold. Let $\T$ be a family of trajectories of (\ref{eq:ss})-(\ref{eq:so}) which is invariant w.r.t. a closed covering $\chi$. Then the pair $(h,f)$ is WZSD w.r.t. $\T$ if each reduced limiting control system for $\T$ is WZSD.
 \end{teo}
\begin{IEEEproof} Suppose that the pair $(h,f)$ is not WZSD w.r.t. $\T$. Then there exists a $\T$-zeroing-output sequence $\{(t_k,x_k,\sigma_k)\}$ of (\ref{eq:ss})-(\ref{eq:so}). Due to Lemma \ref{lm:subseq} there exists a subsequence $\{(x_{k_l},\sigma_{k_l})\}$ which converges to a  limiting trajectory $(\bar x,\bar u)$ of $\T$. Since $\{(t_{k_l},x_{k_l},\sigma_{k_l})\}$ is a $\T$-zeroing-output sequence, $(\bar x,\bar u)$ is a $\T$-zeroing-output trajectory of (\ref{eq:ss})-(\ref{eq:so}). From Theorem \ref{teo:limiting}, there exists a reduced limiting control system $\Sigma_{\gamma}$ for which $(\bar x,\bar u)$ is a bounded complete trajectory. Since $\Sigma_{\gamma}$ is WZSD, it follows that $\inf_{t\in \R}|\bar x(t)|=0$. On the other hand, since $\{(t_k,x_k,\sigma_k)\}$ is a $\T$-zeroing-output sequence of (\ref{eq:ss})-(\ref{eq:so}), there exists $\varepsilon_0>0$ such that $|x_k(t+t_k)|\ge \varepsilon_0$ for all $t\in [-k,k]$ and all $k$. Consequently, for all $t\in \R$, $|\bar x(t)|=\lim_{l\to \infty} |x_{k_l}(t_{k_l}+t)|\ge \varepsilon_0$. Hence we have arrived to a contradiction. Therefore  $(h,f)$ is WZSD w.r.t. $\T$.
\end{IEEEproof}
\section{A characterization for GUAS}\label{sec:guas}
In this section we will give a characterization for the GUAS of a family of trajectories of the switched system (\ref{eq:ss})-(\ref{eq:so}) in terms of the WZSD of its reduced limiting control systems.

By combining Theorems \ref{thm:1} and \ref{thm:criterio_wzsd} it easily follows that if a family of trajectories $\T$ of (\ref{eq:ss})-(\ref{eq:so}) is GUS, satisfies Assumption \ref{ass:integral} and its reduced limiting control systems are WZSD, then it is GUAS. The question which naturally arises is whether the converse of this result holds. We will show that this question has a positive answer if some additional conditions are satisfied.

In what follows we will consider the trivial covering $\chi=\{\chi_i\}_{i=1}^N$, where $\chi_i=\R^n$ for all $i$ and suppose that the family $\T$ satisfies the following.
\begin{as} \label{ass:product} \rm There is a family of switched signals $\tilde \S$ such that $\T$ is the set of all the pairs $(x,\sigma)$ with $\sigma \in \tilde \S$ and $x$ a maximal solution of (\ref{eq:ss}) corresponding to $\sigma$.
\end{as}
We will also assume that $f$ and $h$ in, respectively, (\ref{eq:ss}) and (\ref{eq:so}), satisfy the following conditions, which are slightly stronger than those in Assumption \ref{ass:aap}.
\begin{as} \label{ass:aapstrong}\rm $f$ and $h$ in (\ref{eq:ss})-(\ref{eq:so}) satisfy for all $i\in \I$ the following:
\begin{enumerate}
	\item $h_i(t, \xi)$ is uniformly bounded, measurable in $t$, locally Lipschitz in $\xi$ uniformly in $t$, precompact and $h_i(t,0)\equiv 0$;
		\item $f_i=\hat f_i+\Delta f_i$ with
	\begin{enumerate}
		\item $\hat{f}_i(t,\xi)$ uniformly bounded,  measurable in $t$, locally Lipschitz in $\xi$ uniformly in $t$, and precompact;
		\item $\Delta f_i$ such that for each $R>0$ there exists a constant $M_R>0$ such that $|\Delta f_i(t,\xi)|\le M_R|h_i(t,\xi)|$ for all $t\ge 0$ and all $\xi \in \R^n$ such that $|\xi|\le R$.
	\end{enumerate}	
	\end{enumerate}
\end{as}
\begin{teo} \label{thm:converse} \rm  Let $\T$ be a set of trajectories of (\ref{eq:ss}) for which Assumption \ref{ass:product} holds. Suppose that $f$ in (\ref{eq:ss}) and $h$ in (\ref{eq:so}) satisfy Assumptions \ref{ass:boundedness} and \ref{ass:aapstrong}. Then, the following are equivalent:
\begin{enumerate}
	\item $\T$ is GUAS.
	\item $\T$ is GUS, Assumption \ref{ass:integral} holds and each reduced limiting control system for $\T$ is WZSD.
\end{enumerate}
\end{teo}
\begin{IEEEproof} Since Assumption \ref{ass:aapstrong} implies Assumption \ref{ass:aap}, 2) $\Rightarrow$ 1) straightforwardly follows from Theorems \ref{thm:1} and \ref{thm:criterio_wzsd}.

Next we prove that 1) $\Rightarrow$ 2). Since $\T$ is GUAS, it is GUS. That Assumption \ref{ass:integral} holds can be proved along the lines of the proof of \cite[Prop.4]{leejia_tac05}.

Let $\gamma=\{t_k\}\in \Lambda$ and let $\Sigma_{\gamma}$ be the corresponding reduced limiting control system for $\T$. Let $(\bar x,\bar u)$ be a bounded complete trajectory of $\Sigma_{\gamma}$. Then for all $t\in \R$,
$\bar x(t)=\bar x(0)+\int_0^t \hat F_{\gamma}(s,\bar x(s))\bar u(s)\:ds$
and $\hat H_{\gamma}(s,\bar x(s))\bar u(s)=0$ a.e. on $\R$. Due to Assumption \ref{ass:product}, there exists a sequence $\{\sigma_k\}$ in $\tilde \S$ such that $v_k(\cdot):=\sigma(t_k+\cdot)\rightharpoonup \bar u$. Let $\hat H$ be as in the proof of Theorem \ref{teo:limiting}. For each $k\in \N$ let $g_k(t):=\hat H(t_k+t,\bar x(t))v_k(t)$ for all $t\ge 0$. The functions $g_k$ are nonnegative and uniformly bounded since the components of $\hat H$ are uniformly bounded and $\bar x$ is bounded. Fix $T>0$. Writing
\begin{multline*}
g_k(s)=[\hat H(t_k+s,\bar x(s))-\hat H_{\gamma}(s,\bar x(s))]v_k(s)\\
+ \hat H_{\gamma}(s,\bar x(s))(v_k(s)-\bar u(s))+ \hat H_{\gamma}(s,\bar x(s))\bar u(s),
\end{multline*}
and using the convergence of $\hat H(t_k+\cdot,\bar x(\cdot))$ to $\hat H_{\gamma}(\cdot,\bar x(\cdot))$, the weak convergence of $v_k$ to $\bar u$, the boundedness of the functions $g_k$ and $v_k$ and Lebesgue's convergence theorem, we have that $\lim_{k \to \infty}\int_0^T g_k(s)\:ds =0.$ Since $g_k\ge 0$ for all $k$, it follows that $g_k\to 0$ in $L^1([0,T])$ and therefore, there exists a subsequence of $\{g_k\}$ which converges to $0$ a.e. on $[0,T]$. By using the diagonal process we can extract a subsequence of $\{g_k\}$, which we still denote by $\{g_k\}$, such that $g_k\to 0$ a.e. on $\R_{\ge 0}$. Then, we have that for almost all $t\in \R_{\ge 0}$,
\begin{multline}\label{eq:h-zero}
\lim_{k \to \infty}|h(t_k+s,\bar x(s),v_k(t))|=\\
 \lim_{k \to \infty}\hat H(t_k+s,\bar x(s))v_k(t)=0
\end{multline}
Taking into account the convergence of $\hat F(t_k+\cdot,\bar x(\cdot))$ to $F_{\gamma}(\cdot,\bar x(\cdot))$, the weak convergence of $v_k$ to $\bar u$, the boundedness of $\bar x$ and $\sigma_k$, and the Lebesgue convergence theorem, we have that for every $t\in \R_{\ge 0}$ $\lim_{k \to \infty} \int_0^t \hat F(t_k+s,\bar x(s))v_k(s)\:ds = \int_{0}^{t} \hat F_{\gamma}(s,\bar x(s))\bar u(s)ds.$
Therefore, for all $t\ge 0$,
\begin{align*}
\bar x(t)=\bar x(0)+\lim_{k \to \infty} \int_0^t \hat F(t_k+s,\bar x(s))v_k(s)\:ds.
\end{align*}

Let $\{(x_k,\sigma_k)\}$ be a sequence of trajectories of $\T$ such that $x_k(t_k)=\bar x(0)$ for all $k$ (such a sequence exists due to Assumption \ref{ass:product}). Since $\T$ is GUAS, for every $k$, the trajectory  $(x_k,\sigma_k)$ is forward complete and $|x_k(t)|\le \beta(|\bar x(0)|,t-t_k)$ for all $t\ge t_k$, where $\beta$ is the function of class $\kl$ appearing in \ref{def:guas}) of Definition \ref{def:stab}. In consequence $|x_k(t)|\le R:=\beta(|\bar x(0)|,0)$ for all $t\ge t_k$ and all $k$. Let $z_k(\cdot)=x_k(t_k+\cdot)$ for all $k$. By applying a slight variation of Lemma \ref{def:limsol} we can assume, by passing to a subsequence if necessary, that  there  exists a continuous function $z:\R_{\ge 0}\to \R^n$ such that $z_k \to z$ uniformly on $[0,T]$ for all $T>0$.
Thus, for all $t\ge 0$,
\begin{multline*}
z(t)=\bar x(0)+\lim_{k \to \infty} \int_0^t [\hat F(t_k+s,z_k(s)) \\+\Delta  F(t_k+s,z_k(s))]v_k(s)\:ds.
\end{multline*}
By using the convergence of $z_k$ to $z$, the boundedness of $z_k$ and $v_k$, the local Lipschitzianity and boundedness of the components of $\hat F$, and Lebesgue's convergence theorem, we have that for every $t\ge 0$
\begin{multline*}
\lim_{k \to \infty} \int_0^t \left[\hat F(t_k+s,z_k(s))-
\hat F(t_k+s,z(s))\right ]v_k(s)\,ds \\=0.
\end{multline*}
Therefore, for every $t\ge 0$
\begin{multline}\label{eq:z}
z(t)=\bar x(0)+\lim_{k \to \infty} \int_0^t [\hat F(t_k+s,z(s)) \\+\Delta  F(t_k+s,z_k(s))]v_k(s)\:ds.
\end{multline}
Then, for each $t\ge 0$
\begin{multline}\label{eq:zx}
z(t)-\bar x(t)=\lim_{k\to \infty} \int_0^t [\hat F(t_k+s,z(s))-\hat F(t_k+s,\bar x(s)) \\+\Delta  F(t_k+s,z_k(s))]v_k(s)\:ds.
\end{multline}
Due to the local Lipschitzianity of the functions $\hat f_i$ and the boundedness of $z$ and $\bar x$ there exists a constant $L_1\ge 0$ such that for all $s\ge 0$,
\begin{equation*}
\left|[\hat F(t_k+s,z(s))-\hat F(t_k+s,\bar x(s))]v_k(s)\right |
\le
L_1|z(s)-\bar x(s)|.
\end{equation*}
Then, for all $t\ge 0$,
%%%%%%%%%%%%%%%
\begin{multline}
\limsup_{k\to \infty} \left  | \int_0^t [\hat F(t_k+s,z(s))-\hat F(t_k+s,\bar x(s))]v_k(s)\:ds \right | 
\\ \le
\int_0^t L_1|z(s)-\bar x(s)|\:ds.  \label{eq:L1}
\end{multline}
The boundedness of $z_k$ and $\bar x$ and the local Lipschitzcianity of $h$ imply the existence of $L_2\ge 0$ such that for all $s\ge 0$
\begin{align} 
\zeta_k(s)&=||h(t_k+s,z_k(s),v_k(s))|-|h(t_k+s,\bar x(s),v_k(s))|| \nonumber \\ \label{eq:lips-h}
&\le  L_2|z_k(s)-\bar x(s)|.
\end{align}
From (\ref{eq:lips-h}), the boundedness of $z_k$, 2b) of Assumption \ref{ass:aapstrong}, (\ref{eq:h-zero}), and the convergence of $z_k$ to $z$, it follows that for almost all $s\ge 0$,
\begin{multline}
\limsup_{k\to \infty}|\Delta  F(t_k+s,z_k(s))v_k(s)|=\\ \limsup_{k\to \infty}|\Delta f(t_k+s,z_k(s),v_k(s))|\\ \le
\limsup_{k\to \infty}M_R |h(t_k+s,z_k(s),v_k(s))|
\le \limsup_{k\to \infty} M_R \zeta_k(s) \\
\le \limsup_{k\to \infty} M_R L_2 |z_k(s)-\bar x(s)|=
M_R L_2 |z(s)-\bar x(s)|. \label{eq:ML}
\end{multline}
We then have, from (\ref{eq:ML}) and  Fatou's lemma,
\begin{multline} \label{eq:ML2}
\limsup_{k\to \infty}\int_0^t |\Delta  F(t_k+s,z_k(s))v_k(s)|\:ds\le \\  \int_0^t \limsup_{k\to \infty}|\Delta  F(t_k+s,z_k(s))v_k(s)|ds \\
\le \int_0^t M_RL_2  |z(s)-\bar x(s)|\:ds
\end{multline}
Then, from (\ref{eq:zx}), (\ref{eq:L1}) and (\ref{eq:ML2}), we have, for all $t\ge 0$, that
$$ |z(t)-\bar x(t)|\le (L_1+M_RL_2)\int_0^t|z(s)-\bar x(s)|ds.$$
By applying Gronwall's lemma, we conclude that $z(t)=\bar x(t)$ for all $t\ge 0$. Since $|z(t)|\le \beta(|\bar x(0)|,t)$ for all $t\ge 0$, we have that $z(t)\to 0$ as $t\to \infty$. In consequence $\inf_{t \in \R}|\bar x(t)|=0$ and the WZSD of $\Sigma_{\gamma}$ follows.
\end{IEEEproof}
\begin{remark} \rm Theorem \ref{thm:converse} specialized to the case of nonswitched systems recovers the equivalence 1) $\Leftrightarrow$ 4) in \cite[Thm.3]{leejia_tac05}. \mer 
 \end{remark}
\section{A criterion for UAS and GUAS}
\label{sec:criteria}
In this section we give a criterion for the UAS or GUAS of a family of trajectories of (\ref{eq:ss}). This criterion assumes the existence of weak multiple Lyapunov functions for a family $\T$ which is invariant w.r.t. a closed covering $\chi$. 
 
 \begin{as} \rm \label{ass:wlf} There exists a function $V:\R_{\ge 0}\times \X \times \I \to \R_{\ge 0}$ such that:
\begin{enumerate}
 \item \label{eq:alpha} there exist $\phi_1,\phi_2\in \kk$, such that for all $t\ge 0$, all $\xi\in \chi_i \cap \X$ and all $i\in \I$
 \begin{align}
\phi_1(|\xi|)\le V(t,\xi,i)\le \phi_2(|\xi|);  
 \end{align}
\item  for all $i\in \I$, $V_i(\cdot,\cdot)=V(\cdot,\cdot,i)$ is continuously differentiable and for all $t\ge 0$, all $\xi\in \chi_i \cap \X$ and all $i\in \I$ 
 \begin{align} \label{eq:dotv}
  \dot{V}_i(t,\xi):=\frac{\partial V_i}{\partial t}+ \frac{\partial V_i}{\partial \xi}f_i(t,\xi)\le -\eta_i(t,\xi), 
 \end{align}
 where $\eta_i:\R_{\ge 0}\times \X \to \R_{\ge 0}$ is Lebesgue measurable in $t$ and continuous in $\xi$;
 \item \label{eq:v} for every $(x,\sigma)\in \T$ and every $t_0(x)\le t<s<t_f(x)$ such that $\sigma(t)=\sigma(s)=i$ we have that
 $$V_i(t,x(t))\ge V_i(s,x(s)).$$ 
\end{enumerate} 
\end{as}

The following criterion, which can be seen as an extension of the Krasovskii-LaSalle theorem, is a consequence of Theorems \ref{thm:1} and \ref{thm:criterio_wzsd} and well-known results of the stability theory of switched systems.
\begin{teo}\rm  \label{thm:guas1} Let $\T$ be a family of trajectories of (\ref{eq:ss}) which is invariant w.r.t. a closed covering $\chi$. Suppose that Assumption \ref{ass:wlf} holds and let $h$ in (\ref{eq:so}) be defined by $h(t,\xi,i)=\eta_i(t,\xi)$ for all $t\ge 0$, all $\xi\in \X$ and all $i\in \I$.
Suppose in addition that $f$ in (\ref{eq:ss}) and $h$ satisfy Assumptions \ref{ass:boundedness} and \ref{ass:aap}.

Then $\T$ is UAS if each reduced limiting control system for $\T$ is WZSD. If, in addition, $\X=\R^n$ and the functions $\phi_i$, $i=1,2$, in 1) of Assumption \ref{ass:wlf} are of class $\ki$, then $\T$ is GUAS.
\end{teo}
\begin{remark} \rm As a difference with most UAS results based on multiple weak Lyapunov functions, Theorem \ref{thm:guas1} does not involve any dwell-time assumption on the switching signals of the trajectories of the family $\T$. In fact, it can be used for proving the UAS of families of trajectories whose switching signals do not satisfy any dwell-time condition (see Examples \ref{ex:1} and \ref{ex:4}). 
\mer
\end{remark}
\begin{remark}
	\rm Theorem \ref{thm:guas1} contains as a particular case the generalization of the Krasovskii-LaSalle theorem proposed in \cite[Thm.2]{leejia_tac05} for nonswitched systems. 
\mer
\end{remark}
\begin{IEEEproof}[Proof of Theorem \ref{thm:guas1}] From Proposition 1 and Lemma 6 in \cite{leejia_tac08}, it follows that Assumption \ref{ass:wlf} implies that $\T$ is US, and that it is GUS when $\X=\R^n$ and the functions $\phi_i$, $i=1,2$ in 1) of Assumption \ref{ass:wlf} are of class $\ki$.

Following the steps of the proof of Theorem 2 in \cite{leejia_tac08} it can be proved that there exists an open neighborhood $\mathcal{O}\subset \X$ of the origin   such that $h$ and $\T$ satisfy Assumption \ref{ass:integral} with $\alpha(r)=r$ and $\mathcal{O}$ in place of $\X$, and that we can take $\mathcal{O}=\R^n$ when $\X=\R^n$ and the functions $\phi_1$ and $\phi_2$ in 1) of Assumption \ref{ass:wlf} belong to $\ki$.

Finally, the pair $(h,f)$ is WZSD w.r.t. $\T$ since each reduced limiting control system $\Sigma_{\gamma}$ for $\T$ is assumed WZSD. Therefore, the hypotheses of Theorem \ref{thm:1} are fulfilled (with $\mathcal{O}$ in place of $\X$) and then Theorem \ref{thm:guas1} follows.
\end{IEEEproof}

Next we give a sufficient condition for the uniform asymptotic stability of (\ref{eq:ss}) under arbitrary switchings, that is, when the family of trajectories $\T$ is that of all the possible trajectories of the switched system (\ref{eq:ss}). This condition is a straightforward consequence of Theorem \ref{thm:guas1} and assumes the existence of a common weak Lyapunov function.
\begin{corol} \label{cor:guas_as} \rm
Suppose that $f$ in (\ref{eq:ss}) satisfies Assumption \ref{ass:boundedness} and Assumption \ref{ass:wlf} with a function $V(t,\xi,i)\equiv \bar V(t,\xi)$ and $\chi_i=\R^n$ for all $i\in \I$. Let $h$ in (\ref{eq:so}) be defined by $h(t,\xi,i)=\eta_i(t,\xi)$ for all $t\ge 0$, all $\xi\in \X$ and all $i\in \I$ and suppose in addition  that $f$ and $h$ satisfy Assumption \ref{ass:aap}.

Then (\ref{eq:ss}) is UAS under arbitrary switchings if for each $\gamma \in \Lambda$ the corresponding reduced limiting control system: 
\begin{align} 
\label{eq:lcs_as}
\Sigma_{\gamma}: \;\left \{ \begin{array}{l} \dot{x} = \hat{F }_{\gamma}(t, x) u  \\
y = \hat {H}_{\gamma}(t,x)u=0
\end{array}
\right .,\quad u \in \U,
\end{align}
is WZSD. If, in addition, $\X=\R^n$ and the functions $\phi_i$, $i=1,2$, in 1) of Assumption \ref{ass:wlf} are of class $\ki$, then (\ref{eq:ss}) is GUAS under arbitrary switchings.
\end{corol}
\section{Examples}\label{ej}
In this section we will apply the theory developed in the precedent sections to some interesting examples.    
\begin{ex} \label{ex:1} \rm Consider the switched system (\ref{eq:ss}) in $\X=\R^2$ with $\I=\{1,2,3\}$ and 
 \begin{gather*} 
  f_1(t,\xi)=\left[ \begin{smallmatrix} -g_1(t,\xi)\xi_2 \\
  g_1(t,\xi) \xi_1-\xi_2
                    \end{smallmatrix} \right ],\;
                     f_2(t,\xi)=\left [\begin{smallmatrix} g_2(t,\xi)\xi_2-\xi_1 \\
                      -g_2(t,\xi) \xi_1
                     \end{smallmatrix} \right ] \\ {\rm and} \;
                     f_3(t,\xi)=\left [ \begin{smallmatrix} 2 g_2(t,\xi)\xi_2-\xi_1 \\
                     -2 g_2(t,\xi) \xi_1
                     \end{smallmatrix} \right ].
 \end{gather*}
Suppose that the functions $g_i:\R_{\ge 0}\times \R^2 \to \R$ are uniformly bounded, measurable in $t$ and continuous in $\xi$, uniformly in $t$. Also suppose that for $i=1,2$, $\hat g_i:\R_{\ge 0}\times \R \to \R$, $\hat g_i(t,v)=g_i(t,v e_i)$, where $e_1=[1\;0]'$ and $e_2=[0,1]'$, is precompact and satisfies the following condition: for each $v_0\in \R\setminus\{0\}$, there exist $T=T(v_0)>0$, and $r_i=r_i(v_0)>0$ such that
\begin{align} \label{eq:PE} 
\liminf_{t\to \infty}\int_{t}^{t+T} |\hat g_i(s,v_0)|^{r_i}ds >0.
 \end{align}

We will prove that (\ref{eq:ss}) is GUAS under arbitrary switchings by checking that the switched system satisfies the hypotheses of Corollary \ref{cor:guas_as}. Since the functions $g_i$ are uniformly bounded the same holds for the functions $f_i$, so $f$ satisfies Assumption \ref{ass:boundedness}. It is easy to see that Assumption \ref{ass:wlf} holds with $V(t,\xi,i)\equiv |\xi|^2/2$, $\chi_i=\R^2$, $i\in \I$, the $\ki$-class functions $\phi_j(s)=s^2/2$, $j=1,2$, $\eta_1(t,\xi)=\xi_2^2$ and $\eta_2(t,\xi)=\eta_3(t,\xi)=\xi_1^2$. Let $h_i=\eta_i$ for $i\in \I$. Then Assumption \ref{ass:aap} holds with $\hat f_1(t,\xi)= [0\;\:\hat g_1(t,\xi_1)\xi_1]'$, $\hat f_2(t,\xi)= [\hat g_2(t,\xi_2)\xi_2\;\:0]'$, $\hat f_3(t,\xi)=2\hat f_2(t,\xi)$ and $h$. The functions $\hat f_i$ are precompact because the functions $\hat g_i$ are assumed precompact and the other functions do not depend on $t$.   

Now we proceed to show that each reduced limiting control system of (\ref{eq:ss}) is WZSD. Since $h_i$ does not depend on $t$, $h_{i,\gamma}=h_i$ for all $i\in \I$. On the other hand, from the form of $\hat f_i$ and the fact that $\hat g_i$ is assumed precompact for all $i$, the functions $\hat f_{i,\gamma}$ have the form: $\hat f_{1,\gamma}(t,\xi)= [0\;\hat g_{1,\gamma}(t,\xi_1)\xi_1]'$, $f_{2,\gamma}(t,\xi)= [\hat g_{2,\gamma}(t,\xi_2)\xi_2\;0]'$ and $\hat f_{3,\gamma}(t,\xi)=2\hat f_{2,\gamma}(t,\xi)$, where, for $i=1,2$,  $\hat g_{i,\gamma}$ is the limiting function of $\hat g_i$ associated to $\gamma$. Consequently, for each $\gamma=\{t_k\}\in \Lambda$, the reduced limiting control systems $\Sigma_{\gamma}$ has the form
\begin{align}
 \dot{x}_1 & =(u_2+2u_3)\hat g_{2,\gamma}(t,x_2)x_2 \label{eq:eq1} \\
 \dot{x}_2 & = u_1 \hat g_{1,\gamma}(t,x_1)x_1 \label{eq:eq2}\\
 y & =u_1 x_2^2+(u_2+u_3) x_1^2=0 \label{eq:eq3}
\end{align}
with $u\in \U$. 

Let $(x,u)$ be a bounded complete trajectory of $\Sigma_{\gamma}$. Let $E$ be the set of times $t\in \R$ for which the equalities (\ref{eq:eq1})-(\ref{eq:eq3}) hold, and  let $J_i=\{t\in E:\:u_i(t)>0\}$ for $i\in \I$. Note that $x_2(t)=0$ for all $t\in J_1$ and $x_1(t)=0$ for all $t\in J_2\cup J_3$, and, due to continuity, $x_2(t)=0$ for all $t\in \overline{J_1}$ and $x_1(t)=0$ for all $t\in \overline{J_2\cup J_3}$.

{\em Case 1}. $|J_2\cup J_3|=0$. In this case, $u_1(t)=1$ and $u_2(t)=u_3(t)=0$ a.e. on $\R$, and then $x_2(t)=0$ and $\dot x_1(t)=0$ a.e. on $\R$. Thus $x_1(t)\equiv\theta$ for some $\theta \in \R$ and $0=\dot{x}_2(t)=\hat g_{1,\gamma}(t,\theta)\theta$  a.e. on $\R$. If $\theta \neq 0$, then $\hat g_{1,\gamma}(t,\theta)=0$ for almost all $t\in \R$. Since $\lim_{k \to \infty} \hat g_1(t_k+t,\theta)=\hat g_{1,\gamma}(t,\theta)=0$ for almost all $t\in \R$, and $\hat g_1$ is uniformly bounded, applying the Lebesgue convergence theorem, we have that for $T=T(\theta)$,
\begin{multline*} \liminf_{t\to \infty}\int_{0}^{T} |\hat g_1(t+s,\theta)|^{r_1}ds \\ \le \lim_{k \to \infty}\int_{0}^{T} |\hat g_1(t_k+s,\theta)|^{r_1}ds=0,
\end{multline*}
which contradicts (\ref{eq:PE}). Then $\theta=0$ and $x(t)=0$ for all $t\ge 0$.

{\em Case 2}. $|J_2\cup J_3|>0$ and $|J_1|>0$. Since $|\R\setminus(\cup_{i\in \I}J_i)|=0$, $\overline{J_1}\cup \overline{J_2\cup J_3}=\R$. By the connectedness of $\R$ it follows that $\overline{J_1}\cap \overline{J_2\cup J_3}\neq \emptyset$. In consequence there exists $t^*\in \R$ such that $x_1(t^*)=x_2(t^*)=0$, and then $\inf_{t\in \R}|x(t)|=0$.   

{\em Case 3}. $|J_1|=0$. In this case  $u_1(t)=0$ and $u_2(t)+u_3(t)=1$ for almost all $t$, and then $x_1(t)=0$ and $\dot x_2(t)=0$ for almost all $t$. We then have that $x_2(t)=\theta$ for some constant $\theta$ and $0=\dot{x}_1(t)= [u_2(t)+2u_3(t)]\hat g_{2,\gamma}(t,\theta)\theta$ for almost all $t\in \R$. Since   $u_2(t)+2u_3(t)>0$ for all $t$, $\hat g_{2,\gamma}(t,\theta)\theta=0$ for almost all $t\in \R$. By proceeding as in Case 1, it follows that $\theta=0$ and then, that $x(t)=0$ for all $t\in \R$.

We have then shown that $\Sigma_{\gamma}$ is WZSD for every $\gamma\in \Lambda$. Since the hypotheses of Corollary \ref{cor:guas_as} are fulfilled, we can assert that the switched system is GUAS under arbitrary switchings.
\end{ex}
\begin{remark} \rm The stability of the switched system considered in Example \ref{ex:1}, without mode 3, was studied in \cite[Sec.V]{lee_tan_nesic_tac15}, \cite[Ex.1]{lee_tan_mareels_tac17} and  \cite[Sec.4]{teel_matrosov_automatica16}, assuming stronger regularity conditions on the functions $g_i$ than those we have considered in Example \ref{ex:1} (differentiability in \cite{lee_tan_nesic_tac15, teel_matrosov_automatica16} and almost asymptotic periodicity in \cite{lee_tan_mareels_tac17}). As was asserted in \cite{lee_tan_mareels_tac17}, the fact that there is no common strict quadratic Lyapunov function for the system makes its stability analysis challenging. In \cite{lee_tan_nesic_tac15} the GUAS under arbitrary switchings was proved after performing a series of steps which included the  construction of additional Lyapunov-like functions. In \cite{teel_matrosov_automatica16} the switched system was embedded into a differential inclusion and its GUAS was proved by constructing additional Matrosov functions. To construct Lyapunov-like or Matrosov functions is sometimes a cumbersome and difficult task. In \cite{lee_tan_mareels_tac17} the GUAS of the switched system was proved with the help of a common zeroing output system (CZOS) and a simplified detectability condition. When one considers also mode 3, the resulting switched system has no CZOS and therefore the method developed in \cite{lee_tan_mareels_tac17} cannot be applied. Since the functions $g_i$ are not assumed differentiable, neither the additional Lyapunov function nor the Matrosov functions proposed in \cite{lee_tan_nesic_tac15} and \cite{teel_matrosov_automatica16}, respectively, are differentiable and then they cannot be employed for establishing the GUAS of the switched system as was done in those paper. Example \ref{ex:1} clearly shows the potential of our method as it is able to establish the GUAS under arbitrary switchings of the switched system under weaker conditions and in a straightforward way.        \mer    
 \end{remark}

The following example shows how to exploit the fact that the reduced limiting control systems allow us to take into account the state-dependent constraints the trajectories of the family under analysis satisfy.

\begin{ex}
 \label{ex:4} \rm Consider the switched system (\ref{eq:ss}) in $\R^2$ with three modes, i.e. $\I=\{1,2,3\}$, given by
\begin{gather*} f_1(t,\xi)=\left[ \begin{smallmatrix} b_1(t)\xi_2 \\  -b_1(t)\xi_1 -\alpha_1(t,\xi_2) \end{smallmatrix} \right ], \;  f_3(t,\xi)=\left[ \begin{smallmatrix} -3\xi_1+5 \xi_2 \\ -5\xi_1+3\xi_2 \end{smallmatrix} \right ] \\{\rm and} \; f_2(t,\xi)=\left[ \begin{smallmatrix} -\alpha_2(t,\xi_1) -b_2(t)\xi_2 \\ b_2(t) \xi_1  \end{smallmatrix} \right ]\end{gather*}
where, 
\begin{enumerate}
\item for $j=1,2$, $b_j:\R_{\ge 0}\to \R$ is measurable, bounded and precompact.
\item $b_1$  satisfies the condition
\begin{align}\label{eq:c1}
\liminf_{t\to \infty} \int_{t}^{t+\tau} |b_1(s)|\:ds>0, \quad \forall \tau>0;
\end{align}
\item there exists $T_0>0$ such that 
\begin{align} \label{eq:c2}
\liminf_{t\to \infty} \int_{t}^{t+T_0} |b_2(s)|\:ds>0; 
\end{align} 
\item for $j=1,2$, $\alpha_j:\R_{\ge 0}\times \R \to \R$, $(t,v) \mapsto \alpha_j(t,v)$, is uniformly bounded, measurable in $t$, continuous in $v$, and there exists a continuous and positive definite function $\rho_j:\R\to \R$ such that $\rho_j(v)\le v\alpha_j(t,v)$ for all $v\in \R$.    
\end{enumerate}	

We will show that the family $\T$ of all the forward complete trajectories of (\ref{eq:ss}) which are invariant w.r.t. the covering $\chi=\{\chi_i\}_{i=1}^3$, where $\chi_1=\chi_2=\{\xi \in \R^2:\xi_1\ge 0\}$ and $\chi_3=\{\xi \in \R^2:\xi_1\le 0\}$ is GUAS. Let $V:\R\times \R^2\times \I\to \R$ be defined as $V_1(t,\xi)=V_2(t,\xi)=5\xi_1^2+5\xi_2^2$ and $V_3(t,\xi)=5\xi_1^2-6 \xi_1 \xi_2+5\xi_2^2$. It is clear that $V$ satisfies Assumption \ref{ass:wlf} with $\eta_1(t,\xi)=\rho_1(\xi_2)$, $\eta_2(t,\xi)= \rho_2(\xi_1)$ and $\eta_3(t,\xi)=0$. If we consider the output (\ref{eq:so}) with $h$ defined by $h_i=\eta_i$ for $i\in \I$, we have that $f$ and $h$ satisfy Assumption \ref{ass:aap} with 
\begin{align*} 
  \hat f_1(t,\xi)=\left[ \begin{array}{c} 0 \\  -b_1(t)\xi_1  \end{array} \right ], \; \hat f_2(t,\xi)=\left[ \begin{array}{c}- b_2(t)\xi_2 \\ 0  \end{array} \right ],
\end{align*}
and $\hat f_3=f_3$. For each $\gamma=\{s_k\} \in \Lambda$, the reduced limiting control system $\Sigma_{\gamma}$ has the form
\begin{align*}
\dot{x}_1 &=  -u_2 b_{2,\gamma}(t)x_2-u_3 (3 x_1-5 x_2) \\
                                    \dot{x}_2 &=  -u_1 b_{1,\gamma}(t)x_1-u_3 (5 x_1-3x_2)\\
y &=  u_1 \rho_1(x_2)+u_2 \rho_2(x_1)=0,
\end{align*}
with $u\in \S^*_{\T}$ and $u(t)=U_{x(t)}$ a.e. on $\R$. Here, for $i=1,2$, $b_{i,\gamma}$ is the limiting function of $b_i$ associated to $\gamma$ and 
$U_{\xi}=\{e_3\}$ if $\xi_1<0$, $U_{\xi}={\rm co}\{e_1,e_2\}$ if $\xi_1>0$ and $U_{\xi}={\rm co}\{e_1,e_2,e_3\}$ if $\xi_1=0$.

Let $(x,u)$ be a bounded complete trajectory of $\Sigma_{\gamma}$ such that $\inf_{t\in \R}|x(t)|> 0$. We consider two possible cases.

{\em Case I.} There exists $\bar t$ such that $x(\bar t) \in {\rm int}(\chi_1)=\{\xi:\xi_1>0\}$. In this case, by continuity, $x(t)\in {\rm int}(\chi_1)$ on some interval $[\bar t,T]$, with $T>\bar t$, and hence $u(t)\in {\rm co}\{e_1,e_2\}$ a.e. on $[\bar t,T]$. Then, for almost all $t\in [\bar t,T]$, $\dot{x}_1(t)=-u_2(t)b_{2,\gamma}(t)x_2(t)$,  $\dot{x}_2(t)=-u_1(t) b_{1,\gamma}(t)x_1(t)$,
$u_1(t)\rho_1(x_2(t))+u_2(t) \rho_2(x_1(t))=0$, $u_1(t)+u_2(t)=1$ and $u_i(t)\ge 0$ for $i=1,2$. Since $x_1(t)>0$ for all $t\in [\bar t,T]$, $\rho_2(x_1(t))>0$ for all $t\in [\bar t,T]$, and then $u_2(t)=0$ a.e. on $[\bar t,T]$. Therefore, for almost all $t\in [\bar t,T]$ we have that $\rho_1(x_2(t))=0$, and in consequence $x_2(t)=0$ for all $t\in  [\bar t,T]$ due to continuity. Then for almost all $t\in  [\bar t,T]$, $0=\dot{x}_2(t)= -b_{1,\gamma}(t)x_1(t)$. Let $Z=\{t\in [\bar t,T]:\:b_{1,\gamma}(t)\neq 0\}$. Since, due to the boundedness of $b_1$, Lebesgue's convergence theorem and (\ref{eq:c1}) we have that
\begin{multline*}\int_{\bar{t}}^{T} |b_{1,\gamma}(s)|\:ds = \lim_{k\to \infty}\int_{\bar{t}}^{T}|b_1(s_k+s)|\:ds \\ \ge \liminf_{\tau\to \infty} \int_{\tau}^{\tau+T-\bar t}|b_1(s)|\:ds>0,\end{multline*}
it follows that $|Z|>0$. Consequently there exists $t^*\in Z$ such that $b_{1,\gamma}(t^*)x_1(t^*)=0$. Therefore $x_1(t^*)=0$ and then $x(t^*)=0$. So case I is impossible.

{\em Case II.} $x(t)\in \chi_3$ for all $t\in \R$. Since while $x(t)$ remains in ${\rm int}(\chi_3)=\{\xi:\xi_1<0\}$, $x$ is a solution of the linear differential equation $\dot{x}=f_3(t,x)$, and since the orbits of that system are clockwise ellipses, there exists $\bar t\in \R$ such that $x(\bar t)\in \partial\chi_3=\{\xi:\xi_1=0\}$. Suppose that $x(t)\in \partial\chi_3$ for all $t\ge \bar t$. Since $x_2(t)\neq 0$ for all $t\ge \bar t$ because we suppose that $x(t)\neq 0$ for all $t\in \R$, it follows that $u_1(t)=0$ for almost all $t\in [\bar t,\infty)$. Since $x_1(t)=0$ for all $t\in [\bar t,\infty)$, it follows that $0=\dot{x}_1(t)= [-u_2(t)b_{2,\gamma}(t)+5 u_3(t)]x_2(t)$ a.e. on $[\bar t,\infty)$. The latter is only possible if $b_{2,\gamma}(t)\ge 0$ a.e. on $[\bar t, \infty)$. In such a case $u_2(t)=5/(5+b_{2,\gamma}(t))$ and $u_3(t)=b_{2,\gamma}(t)/(5+b_{2,\gamma}(t))$ a.e. on $[\bar t,\infty)$. Since $b_{2,\gamma}(t)\le m_2$ for some $m_2\ge 0$, because $b_2$ is bounded, there exists $c>0$ such that $u_3(t)\ge c b_{2,\gamma}(t)$ a.e. on $[\bar t,\infty)$. Taking into account that $\dot{x}_2(t)=3 u_3(t) x_2(t)$ a.e. on $[\bar t,\infty)$, we have that 
$$x_2(t)=x_2(\bar t)e^{3\int_{\bar t}^{t}u_3(s) \:ds}\ge x_2(\bar t)e^{3 c\int_{\bar t}^{t} b_{2,\gamma}(s) \:ds}.$$
Since the integral $\int_{\bar t}^{\infty}b_{2,\gamma}(s) \:ds$ is divergent due to (\ref{eq:c2}) and the fact that $b_{2,\gamma}(t)\ge 0$ for almost all $t\ge \bar t$, it follows that $|x(t)|\to \infty$ as $t\to \infty$, which contradicts the fact that $x$ is bounded. 

Therefore $x$ does not remain in $\partial\chi_3$ for all $t\ge \bar t$. We will consider two cases.\\
i) $x_2(\bar t)>0$. Then there exists $t'\ge \bar t$ such that $x(t) \in \partial \chi_3$ for all $t\in [\bar t,t']$ and a sequence $\{t_k\}$ such that $t_k \searrow t'$ and $x(t_k)\in {\rm int}(\chi_3)$. Note that by continuity and the fact that $x(t)\neq 0$ for all $t\in \R$, $x_2(t')>0$. Since $t_k \searrow t'$ and $x(t_k)\in {\rm int}(\chi_3)$, for each $k$ there exists $\delta_k>0$, with $\delta_k\to 0$ such that $x(t_k-\delta_k)\in \partial \chi_3$ and $x(t)\in {\rm int}(\chi_3)$ for all $t\in (t_k - \delta_k,t_k]$. Then $x_2(t_k-\delta_k)< 0$ for all $k$, since $x$ is a solution of $\dot x=f_3(t,x)$ on $ (t_k - \delta_k,t_k]$ and the nontrivial solutions of this equation are clockwise ellipses. So $x_2(t')=\lim_{k \to \infty}x_2(t_k-\delta_k) \le 0$, which is a contradiction. So, case i) is impossible.

ii) $x_2(\bar t)<0$. Since $x$ leaves $\partial \chi_3$, there exists $t'>\bar t$ such that $x(t')\in {\rm int}(\chi_3)$. Then, by using the facts that while $x$ remains in ${\rm int}(\chi_3)$, $x$ is a solution of the linear differential equation $\dot{x}=f_3(t,x)$, and the solutions of that equation are clockwise ellipses, there exists $t''>t'$ such that $x(t'')\in  \partial \chi_3$ and $x_2(t'')>0$. Then, reasoning as in case i), we conclude that case ii) is also impossible.

We  have then proved that for every bounded complete trajectory $(x,u)$ of the reduced limiting control system $\Sigma_2$ it holds that $\inf_{t \in \R}|x(t)|=0$. So $\Sigma_{\gamma}$ is WZSD for all $\gamma$ and the family $\T$ is GUAS. 
\end{ex}

\begin{remark} \rm To establish the GUAS of the family of trajectories $\T$ considered in Example \ref{ex:4} is challenging. On one hand, it seems difficult to find a family of strict Lyapunov functions $V_i$ that satisfy the conditions in Assumption \ref{ass:wlf}, which are standard conditions for establishing GUAS by means of multiple strict Lyapunov functions. On the other hand, the existing results based on multiple weak Lyapunov functions, namely extensions of LaSalle's invariance principle or of the Krasovskii-LaSalle theorem, cannot be applied to this example, even in the case of time-invariant subsystems, since they assume the trajectories of switched system satisfy some kind of dwell-time constraint, and the trajectories of $\T$ do not satisfy any of them. Neither the results in \cite{lee_tan_mareels_tac17}, which do not make any dwell-time assumption,  can be applied to this example since the system does not have a CZOS. The only results we are aware of, that could be applied to our example are those in \cite{lee_tan_nesic_tac15}, but they require to find an additional Lyapunov-like function $W$ verifying condition (H2) in \cite{lee_tan_nesic_tac15} and such that the pair $((h,\dot W), f)$ be OPE (or WZSD) w.r.t. $\T$, that in this case seems a nontrivial task due to the facts that $f_3$ is stable, but not asymptotically, and the functions $b_i$ neither have a definite sign nor are differentiable. \mer
\end{remark} 

The following example illustrates how to deal with time-dependent constraints on the switching signals of a family of trajectories of a switched system, when analyzing its WZSD via the reduced limiting control systems.
\begin{ex} \label{ex:3} \rm
 Consider the ideal switched model of the semi-quasi-Z-source inverter \cite{caojia_tpe11,denhai_auto16}, connected to a nonlinear time-varying resistive load and under zero input voltage:
  \begin{align}
  \label{eq:ex1}  \dot x &= f(t,x,\sigma) = A_{\sigma} x - e_4 g_{\sigma}(t,e_4'x),
  \end{align}
  where $\sigma:\R\to \{1,2\}$ is a switching signal, $e_4 = [0\ 0\ 0\ 1]'$, $P = {\rm diag}(L_1,L_2,C_1,C_2)$,
  \begin{align*}
    A_1 = P^{-1}\left[ 
      \begin{smallmatrix}
        0 & 0 & 0 & 0\\
        0 & 0 & 1 & 1\\
        0 & -1 & 0 & 0\\
        0 & -1 & 0 & 0
      \end{smallmatrix}\right],\quad
    A_2 = P^{-1}\left[
      \begin{smallmatrix}
        0 & 0 & -1 & 0\\
        0 & 0 &  0 & 1\\
        1 & 0 &  0 & 0\\
        0 & -1 & 0 & 0
      \end{smallmatrix}\right],
    \end{align*}  
 and, for $i=1,2$, $g_i:\R_{\ge 0}\times \R \to \R$, $(t,v)\mapsto g_i(t,v)$, is uniformly bounded, measurable in $t$, continuous in $v$, and there exists  a continuous and positive definite function $\ell_i:\R\to \R$ such that $\ell_i(v)\le v g_i(t,v)$ for all $t\ge 0$ and all $v\in \R$. The positive constants $L_1,L_2,C_1,C_2$ represent the inverter inductance and capacitance values. 
  
  We will show that, irrespective of the load functions $g_i$, the GUAS of this inverter model 
  can be ensured if the switching signals $\sigma$ belong to the set of switching signals $\S[T,d_m,d_M]$, with $T>0$ and $0 < d_m < d_M < \pi \sqrt{L_1C_1}$. We define the set $\S[T,d_m,d_M]$ as follows: $\sigma \in \S[T,d_m,d_M]$ if for each $t\ge 0$ there exist four times $t\le \tau_1^{\sigma}(t)<\tau_2^{\sigma}(t)<\tau_3^{\sigma}(t)<\tau_4^{\sigma}(t)\le t+T$ such that 
  \begin{enumerate}
   \item $d_m\le \tau_{i+1}^{\sigma}(t)-\tau_i^{\sigma}(t)\le d_M$ for $i=1,2,3$;
   \item $\sigma(s)=1$ for all $s\in [\tau_1^{\sigma}(t),\tau_2^{\sigma}(t))\cup [\tau_3^{\sigma}(t),\tau_4^{\sigma}(t))$ and $\sigma(s)=2$ for all $s\in [\tau_2^{\sigma}(t),\tau_3^{\sigma}(t))$.
  \end{enumerate}
Let $\T$ be the set of all the trajectories $(x,\sigma)$ of (\ref{eq:ex1}) such that $\sigma \in \S[T,d_m,d_M]$.
  
It is clear that $f$ verifies Assumption \ref{ass:boundedness}. Let $V(t,\xi,i)= \frac{1}{2}\xi'P\xi$. Then $V$ satisfies 1) of Assumption \ref{ass:wlf} with the $\ki$-class functions $\phi_1(s) = \lambda_{m}s^2$ and $\phi_2(s)=\lambda_{M}s^2$, where $\lambda_{m},\lambda_{M}$ are the minimum and maximum eigenvalues of $P/2$. Item 2) of Assumption \ref{ass:wlf} holds with $\eta_i(t,\xi)=C_2 \ell_i(e_4'\xi)$, because
  \begin{align*}
    \dot{V}_i(t,\xi) =
    \xi'P A_i \xi - \xi' C_2 e_4  g_i(t,e_4'\xi) \le  -C_2 \ell_i(e_4'\xi)\le 0.
    \end{align*}
Since $V$ is a common weak Lyapunv function for the switched system, then 3) of Assumption \ref{ass:wlf} also holds. We then conclude that $\T$ satisfies Assumption \ref{ass:wlf} with the functions $V$, $\phi_1$, $\phi_2$, $\eta_1$ and $\eta_2$. Let, for $i=1,2$, $h_i=\eta_i$. Then $f$ and $h$ satisfy Assumption \ref{ass:aap} with $\hat{f_i}(t,\xi)=\hat A_i \xi$, where
\begin{align*}
   \hat A_1 = P^{-1}\left[ 
      \begin{smallmatrix}
        0 & 0 & 0 & 0\\
        0 & 0 & 1 & 0\\
        0 & -1 & 0 & 0\\
        0 & -1 & 0 & 0
      \end{smallmatrix}\right],\quad
  \hat  A_2 = P^{-1}\left[
      \begin{smallmatrix}
        0 & 0 & -1 & 0\\
        0 & 0 &  0 & 0\\
        1 & 0 &  0 & 0\\
        0 & -1 & 0 & 0
      \end{smallmatrix}\right],
    \end{align*}   
Since the functions $\hat f_i$ and $h_i$ do not depend on $t$, there is only one reduced limiting control system $\Sigma$ which is given by
\begin{align*}
 \Sigma: \left \{ \begin{array}{l} \dot{x} = u_1 \hat A_1 x+ u_2 \hat A_2 x  \\
y = u_1 h_1(t,x)+u_2 h_2(t,x)=0
\end{array}
\right .,\quad u \in \S^*_{\T}.
\end{align*}
where $\S^*_{\T}$ is the set of controls $u\in \U$ for which there exist a sequence $\{t_k\}$ with $t_k\to \infty$ and a sequence $\{\sigma_k\}$ in $\S[T,d_m,d_M]$ such that $\sigma_k(t_k+\cdot) \rightharpoonup u$. Next we will prove that $\Sigma$ is WZSD.

Let $(x,u)$ be a bounded complete trajectory of $\Sigma$. The following claim, whose proof is given in the Appendix, will be used in the following. 

\begin{claim} 
\label{cl:claim} \rm There exist four times $0\le \tau_1<\tau_2<\tau_3<\tau_4\le T$ such that
\begin{enumerate}
   \item $d_m\le \tau_{i+1}-\tau_i\le d_M$ for $i=1,2,3$;
   \item $u(s)=e_1$ for almost all $s\in [\tau_1,\tau_2)\cup [\tau_3,\tau_4)$ and $u(s)=e_2$ for almost all $s\in [\tau_2,\tau_3)$.
  \end{enumerate}
\end{claim}
From Claim \ref{cl:claim} and the fact that $(x,u)$ is a trajectory of $\Sigma$, we have that $x$ satisfies the following:
\begin{itemize}
 \item[i)] $\dot{x}(t)=\hat A_1 x(t)$ and $h_1(t,x(t))=0$ for almost all $t\in [\tau_1,\tau_2)\cup [\tau_3,\tau_4)$;
 \item[ii)] $\dot{x}(t)=\hat A_2 x(t)$ and $h_2(t,x(t))=0$ for almost all $t\in [\tau_2,\tau_3)$;
\end{itemize}
From i), the continuity of $x$ and taking into account that $h_1(t,x(t))=0$ for almost all $t\in [\tau_1,\tau_2)\cup [\tau_3,\tau_4)$ implies that $x_4(t)\equiv 0$ on $[\tau_1,\tau_2)\cup [\tau_3,\tau_4)$,  we have that for all $t\in [\tau_1,\tau_2)\cup [\tau_3,\tau_4)$
\begin{gather*} 
\dot x_1(t)=0,\qquad \dot x_2(t)=\frac{1}{L_2}x_3(t),\qquad \dot x_3(t)=-\frac{1}{C_1}x_2(t)\\ {\rm and} \qquad 0=\dot x_4(t)=-\frac{1}{C_2}x_2(t).
\end{gather*}
From the equations above and using the continuity of $x$ it follows that $x_2(\tau_i)=x_3(\tau_i)=x_4(\tau_i)=0$ for $i=2,3$. Taking into account ii), the continuity of $x$ and the fact that $h_2(t,x(t))=0$ for almost all $t\in [\tau_2,\tau_3)$ implies that $x_4(t)\equiv 0$ on $[\tau_2,\tau_3)$,  we have that for all $t\in [\tau_2,\tau_3)$
\begin{gather*} 
\dot x_1(t)=-\frac{1}{L_1}x_3(t),\qquad \dot x_2(t)=0,\qquad \dot x_3(t)=\frac{1}{C_1}x_1(t) \\ {\rm and} \qquad 0=\dot x_4(t)=-\frac{1}{C_2}x_2(t).
\end{gather*} 
In particular, $x_1$ is a solution of the following second order boundary value problem: 
\begin{gather*}
 \ddot x_1(t)+\frac{1}{L_1C_1} x_1(t)=0, \quad \tau_2<t<\tau_3, \\ \dot x_1(\tau_2)=\dot x_1(\tau_3)=0.
\end{gather*}
Since $\tau_3-\tau_2<\pi \sqrt{L_1C_1}$, it follows that $x_1(t)\equiv 0$ on $[\tau_2,\tau_3]$. In consequence, $x_1(\tau_2)=0$ and then $x(\tau_2)=0$. Consequently $\inf_{t\in \R}|x(t)|=0$. Therefore, $\Sigma$ is WZSD and the GUAS of $\T$ follows from Theorem \ref{thm:criterio_wzsd}, since the hypotheses of  that theorem are fulfilled. 
\end{ex}
\begin{remark} To prove the GUAS of the semi-quasi-Z-source inverter considered in Example \ref{ex:3} is not a trivial task, since although each component subsystem is GUS, neither of them is asymptotically stable and consequently the system is not GUAS for arbitrary switchings. The GUAS of the inverter is obtained via adequate switching as the results in the papers \cite{denhai_auto16} and \cite{manhai_tac17} show. In these papers the GUAS was proved assuming stronger conditions on the switching signals, namely the existence of minimum and maximum dwell-times. In both papers an auxiliary stability result for the inverter (essentially Lemma 1 in \cite{haimid_cdc13}) was used, whose proof was obtained via an {\em ad-hoc} method.    \mer
 \end{remark}
\section{Conclusion} \label{sec:conclusion}
In this paper we studied the (G)UAS of families of trajectories of switched NLTV systems, whose switchings verify time/state dependent constraints, by means of the detectability of output-maps. With this aim, we introduced the notion of reduced limiting control systems for switched NLTV systems and the associated notion of WZSD for those limiting systems. We proved that if the system is (G)US, the output-map satisfies a certain integrability condition and each reduced limiting system is WZSD, then the family of trajectories under study is (G)UAS. For certain classes of families of trajectories, and under slightly stronger hypotheses, these sufficient conditions are also necessary. The obtained sufficient conditions in conjunction with the existence of weak multiple Lyapunov functions allowed us to obtain a (G)UAS criterion for families of trajectories of the system. This criterion can be seen as an extension of the classical Krasovskii-LaSalle theorem. An interesting feature of our approach is that no dwell-time assumptions are needed. In fact, when a common weak Lyapunov function exists, the extension of the Krasovskii-LaSalle theorem enabled us to obtain a criterion for the GUAS under arbitrary switchings. We illustrated the effectiveness of our results by means of several interesting examples, including the stability analysis of a semi-quasi-Z-source inverter.
\section{Appendix}
\subsection{Proof of Lemma \ref{lem:0-der}}
Since the Lebesgue measure is regular, for each $n\in \N$ there exists a closed set $E_n\subset E$ such that $|E\setminus E_n|< 1/n$. We can suppose, due to Exercise 27 in Chapter 2 of \cite{pma_rudin}, that $E_n$ is perfect for every $n$. Let $E^*=\cup_{n\ge 1}E_n$. Then $E^*\subset E$ and $|E\setminus E^*|=0$. Let $t\in E^*$. Then $t\in E_n$ for some $n$. Since $E_n$ is perfect, there exists a sequence $\{t_k\}$ in $E_n$ such that $t_k\to t$ and $t_k\neq t$ for all $k$. Note that $t\in E$ and $t_k\in E$ for all $k$. In consequence $\dot \varphi(t)=\lim_{k\to \infty}(\varphi(t_k)-\varphi(t))/(t_k-t)=0$ since $\varphi(t_k)=0$ for all $k$.
\qed
  
\subsection{Proof of Lemma \ref{lem:sdc}}
To prove Lemma \ref{lem:sdc} we need the following result about the set of control values $U_{\xi}$.
\begin{lema}
	\label{lem:usc} \rm Let $\chi=\{\chi_i\}_{i=1}^N$ be a closed covering of $\R^n$. Then for each $\xi\in \R^n$ there exists $\delta>0$ such that
	$$ U_{\zeta}\subset U_{\xi}\quad \forall \zeta:\;|\zeta-\xi|<\delta.$$
	\end{lema}
\begin{IEEEproof} Let $\I_{\xi}^c=\I\setminus \I_{\xi}$. Then, for each $i\in \I_{\xi}^c$, $\xi \notin \chi_i$ and therefore $\delta=\min_{i\in \I_{\xi}^c}d(\xi,\chi_i)>0$. Here $d(\xi,\chi_i)$ denotes the distance from $\xi$ to $\chi_i$. Let $\zeta\in \R^n$ be such that $|\zeta-\xi|<\delta$. Then $\zeta\notin \chi_i$ for all $i\in \I_{\xi}^c$ and, consequently, $\I_{\zeta}\subset \I_{\xi}$. Thus, $U_{\zeta}\subset U_{\xi}$. 
\end{IEEEproof}

\begin{IEEEproof}[Proof of Lemma \ref{lem:sdc}]
	 Let $(\bar x,\bar u)$ be a limiting trajectory of $\T$. Then there exist a sequence $\{t_k\}$ in $\R_{\ge 0}$ and a sequence $\{(x_k,\sigma_k)\}$ in the conditions of Definition \ref{def:limsol}. Let $z_k$ and $v_k$ as in the proof of Theorem \ref{teo:limiting}. Then $z_k\to \bar x$ uniformly on compact subsets of $\R$ and $v_k \rightharpoonup \bar u$. Since $\T$ is invariant w.r.t. $\chi$, then $v_k(t)\in U_{z_k(t)}$ for all $t\in \R$. Let $T>0$ be arbitrary. Since $v_k \in \U\subset L^{\infty}(\R)$, the restriction of $v_k$ to $[-T,T]$ belongs to $L^1([-T,T])$. Given that  $L^{\infty}([-T,T])\subset L^1([-T,T])$,  $v_k$ converges to $\bar u$ weakly in $L^1([-T,T])$. (Here weak convergence in $L^1([-T,T])$ is the convergence corresponding to the weak topology in $L^1([-T,T])$ induced by its dual $L^{\infty}([-T,T])$). Then, from \cite[Thm. 3.13]{rudin_functional_analysis} it can be derived the existence of a sequence $\{w_k\}$, with $w_k=\sum_{j=k}^{m^k}\lambda_j^k v_j$, $\lambda_j^k\ge 0$ and $\sum_{j=k}^{m_k}\lambda_j^k=1$, such that $w_k\to \bar u$ in the norm of $L^1([-T,T])$. Hence, by passing to a subsequence and relabeling, we can assume that $w_k(t)\to \bar u(t)$ for almost all $t\in [-T,T]$. Let $t\in [-T,T]$ be such that $w_k(t)\to \bar u(t)$. As $v_j(t)\in U_{z_j(t)}$ for all $j$ and since $z_k(t)\to \bar x(t)$, by applying Lemma \ref{lem:usc}, we have that $v_j(t)\in U_{\bar x(t)}$ for $j$ large enough and then $w_k(t)\in U_{\bar x(t)}$ for $k$ large enough. Hence, $\bar u(t)\in U_{\bar x(t)}$ since $U_{\bar x(t)}$ is closed. We have then proved that for all $T>0$, $\bar u(t)\in  U_{\bar x(t)}$ for almost all $t\in [-T,T]$. The fact that $\bar u(t)\in  U_{\bar x(t)}$ for almost all $t\in \R$ easily follows by taking into account that $\R=\cup_{j=1}^\infty[-j,j]$.
 \end{IEEEproof}
\subsection{Proof of Claim \ref{cl:claim}}
In order to prove Claim \ref{cl:claim} we will need the following fact.
\begin{lema} \label{lem:paste} 
 \rm Let $\{u_k\}$ be a sequence in $\U$ such that $u_k\rightharpoonup u$. Suppose there exists a measurable function $v:E\to U$, with $E$ a Lebesgue measurable set, such that $\lim_{k\to \infty}u_k(t)=v(t)$ for almost all $t\in E$. Then $u(t)=v(t)$ for almost all $t\in E$.
 \end{lema}
\begin{IEEEproof} Let $\hat u:\R\to U$ be such that $\hat u=u$ on $\R\setminus E$ and $\hat u=v$ on $E$. Let $f$ be any function in $L^1_N(\R)$. Then
\begin{multline*} \int_{-\infty}^{\infty} f'(t) u_k(t)\:dt=\int_{-\infty}^{\infty} I_{\R\setminus E}(t) f'(t) u_k(t)\:dt \\ + \int_{-\infty}^{\infty} I_{E}(t) f'(t) u_k(t)\:dt,
\end{multline*}
where $I_A$ denotes the indicator function of the set $A$.

Since $u_k \rightharpoonup u$, we have that 
$$ \lim_{k\to \infty} \int_{-\infty}^{\infty} I_{\R\setminus E}(t) f'(t) u_k(t)\:dt=\int_{-\infty}^{\infty} I_{\R\setminus E}(t) f'(t) u(t)\:dt.$$
The pointwise convergence of $u_k$ to $v$ on $E$, the uniform boundedness of $\{u_k\}$ in $L_N^{\infty}(\R)$, the fact that $I_{E}\:f\in L^{1}_N(\R)$ and the Lebesgue convergence theorem imply that
$$ \lim_{k\to \infty} \int_{-\infty}^{\infty} I_{E}(t) f'(t) u_k(t)\:dt=\int_{-\infty}^{\infty} I_{E}(t) f'(t) v(t)\:dt.$$
Therefore
\begin{multline*} 
\lim_{k\to \infty} \int_{-\infty}^{\infty} f'(t) u_k(t)\:dt \\=\int_{-\infty}^{\infty} f'(t)[I_{\R\setminus E}(t)u(t)+I_{E}(t)v(t)]\:dt \\=\int_{-\infty}^{\infty} f'(t) \hat u(t)\:dt.
\end{multline*}
Consequently $u_k \rightharpoonup \hat u$ and then $u=\hat u$ a.e. on $\R$.
\end{IEEEproof}

\begin{IEEEproof}[Proof of Claim \ref{cl:claim}] Since $u\in \S^*_{\T}$, there exist a sequence $\{t_k\}$ with $t_k\to \infty$ and a sequence $\{\sigma_k\}$ in $\S[T,d_m,d_M]$ such that $\sigma_k(t_k+\cdot) \rightharpoonup u$. Since $\{\sigma_k\} \in \S[T,d_m,d_M]$ there exist, for each $k$, four times
$t_k\le \tau_1^{\sigma_k}(t_k)<\tau_2^{\sigma_k}(t)<\tau_3^{\sigma_k}(t)<\tau_4^{\sigma}(t_k)\le t_k+T$ such that 
  \begin{enumerate}
   \item $d_m\le \tau_{i+1}^{\sigma_k}(t_k)-\tau_i^{\sigma_k}(t_k)\le d_M$ for $i=1,2,3$;
   \item $\sigma_k(s)=1$ for all $s\in [\tau_1^{\sigma_k}(t_k),\tau_2^{\sigma_k}(t_k))\cup [\tau_3^{\sigma_k}(t_k),\tau_4^{\sigma_k}(t_k))$ and $\sigma_k(s)=2$ for all $s\in [\tau_2^{\sigma_k}(t_k),\tau_3^{\sigma_k}(t_k))$.
  \end{enumerate}
Let, for all $k$ and $i=1,\ldots,4$, $\tau_i^k=\tau_i^{\sigma_k}(t_k)-t_k$. Since the sequence $\{\tau_i^k\}$ is bounded for all $i$, there exist a strictly increasing subsequence $\{k_j\}$ of $\{k\}$, and real numbers $\tau_i$, $i=1,\ldots, 4$, such that $\tau_i^{k_j}\to \tau_i$ for all $i$. In addition, $0 \le \tau_1<\tau_2<\tau_3<\tau_4\le T$ and $d_m\le \tau_{i+1}-\tau_i\le d_M$. Let $t\in (\tau_1,\tau_2)\cup(\tau_3,\tau_4)$. Then $t_{k_j}+t \in \left(\tau_1^{\sigma_{k_j}}(t_{k_j}),\tau_2^{\sigma_{k_j}}(t_{k_j})\right)\cup \left(\tau_3^{\sigma_{k_j}}(t_{k_j}),\tau_4^{\sigma_{k_j}}(t_{k_j})\right)$ for $j$ large enough. Then $\sigma_{k_j}(t)=e_1$ for $j$ large enough and $\lim_{j\to \infty}\sigma_{k_j}(t_k+t)=e_1$ (here we have interpreted $\sigma_{k_j}$ as the control $u_{\sigma_{k_j}}$). Similarly, we can prove that $\lim_{j\to \infty}\sigma_{k_j}(t_k+t)=e_2$ for all $t\in (\tau_2,\tau_3)$. Since $\sigma_{k_j}(t_{k_j}+\cdot)\rightharpoonup u$, from Lemma \ref{lem:paste}, we have that $u(t)=e_1$ for almost all $t\in  [\tau_1,\tau_2)\cup[\tau_3,\tau_4)$ and $u(t)=e_2$ for almost all $t\in  [\tau_2,\tau_3)$.
\end{IEEEproof}
\balance
\bibliographystyle{IEEEtran}
\bibliography{Basica_bib}
\begin{IEEEbiography}[{\includegraphics[width=1in,height=!]{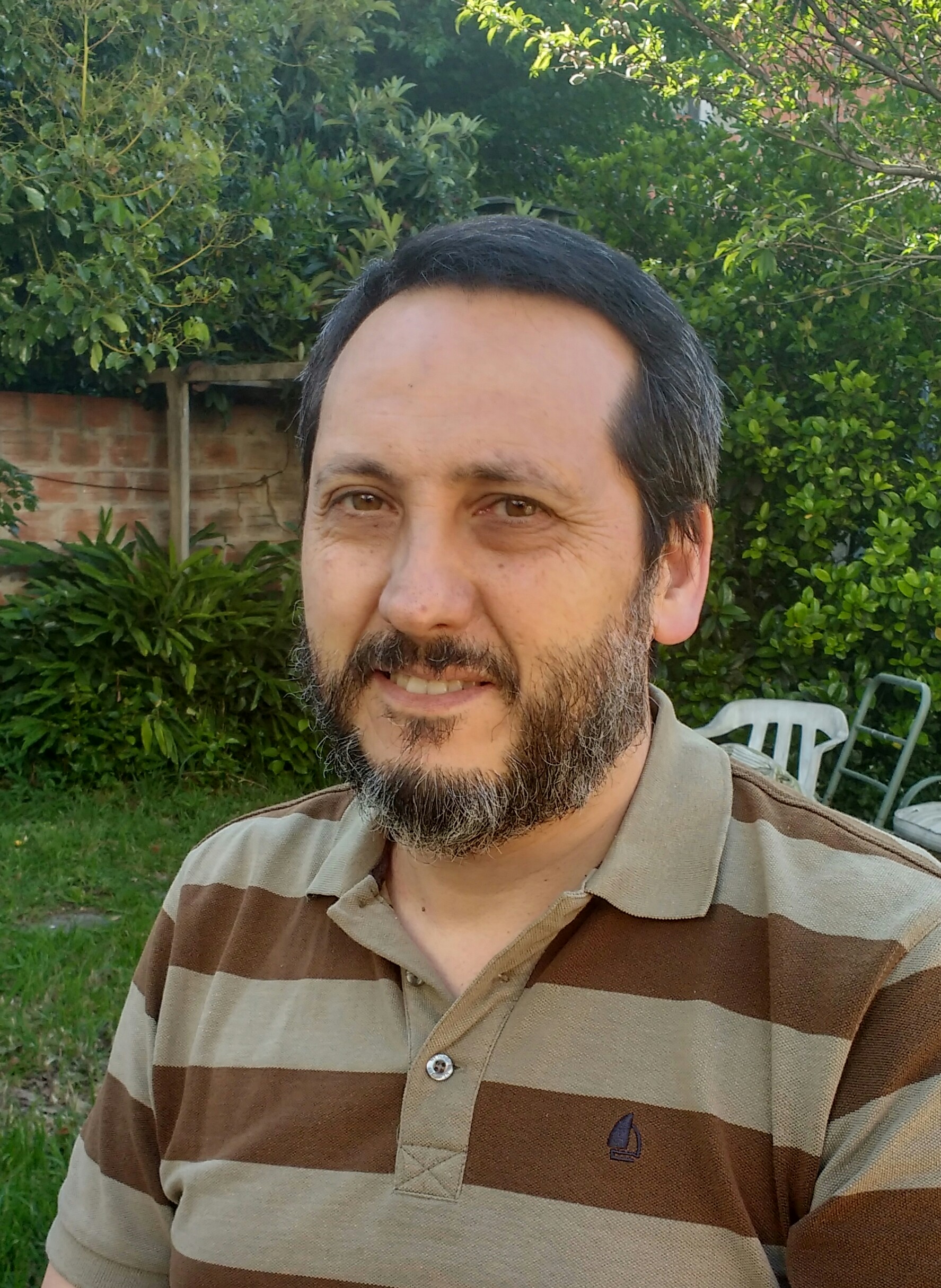}}]
	{Jos\'e Luis Mancilla Aguilar} received the Licenciado en
	Matem\'atica degree (1994) and his Doctor’s degree in
	Mathematics (2001) from the Universidad Nacional de
	Buenos Aires (UBA), Argentina. From 1993 to 1995, he
	received a Research Fellowship from the Argentine Atomic
	Energy Commission (CNEA) in  
	nonlinear control. Since 1995, he has been with the Department of
	Mathematics of the Facultad de Ingenier\'{\i}a (UBA), where he is currently a
	part-time Associate Professor. Since 2005, Dr. Mancilla-Aguilar has held a
	Professor position at the Department of Mathematics of the Instituto
	Tecnol\'ogico de Buenos Aires (ITBA) and currently is the head of the Centro
	de Sistemas y Control (CeSyC). His research interests include hybrid systems
	and nonlinear control.
\end{IEEEbiography}

\begin{IEEEbiography}[{\includegraphics[width=1in,height=!]{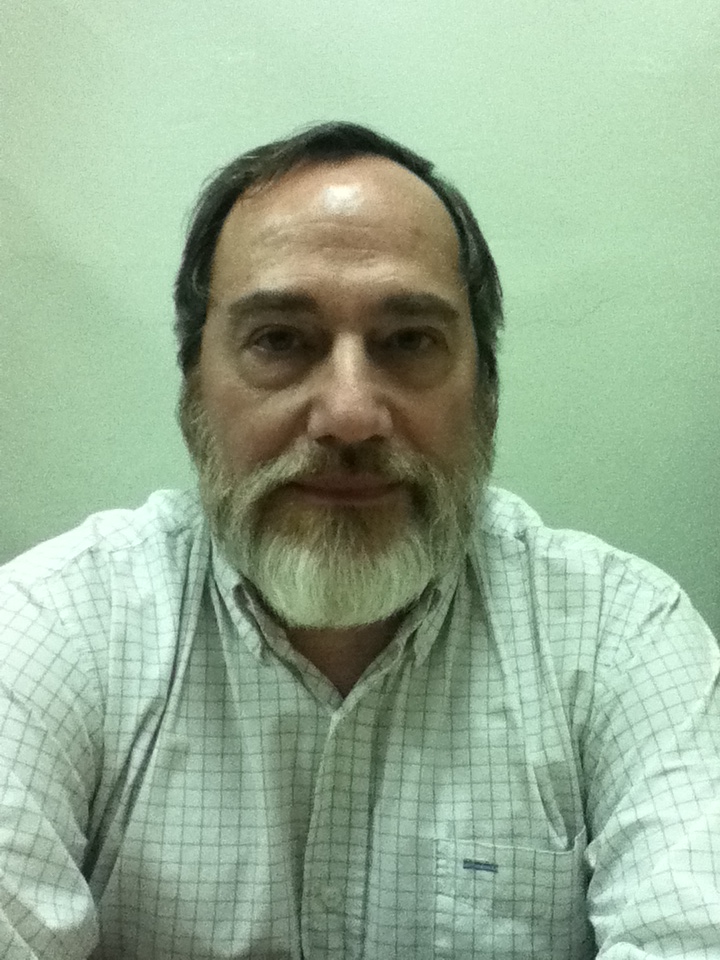}}]
	{Rafael Antonio Garc\'{\i}a} received the Engineering degree in Electronics in 1979, the Licenciado degree in Mathematics in 1984 and the Ph.D. degree, also in Mathematics in 1993, all from the University of Buenos Aires. From 1979 to 1987, he worked in the Instituto de Investigaciones Cient\'{\i}ficas y T\'ecnicas de las Fuerzas Armadas in 
	advanced communications systems. Since 1995 he has been Professor of Mathematics and of Control Theory at the Faculty of Engineering of the University of Buenos Aires, where he is currently a part-time Associate Professor. Since 2002 Dr. Garc\'{\i}a has been the head of the Department of Mathematics of the Instituto Tecnol\'ogico de Buenos Aires (ITBA). His main research interests are in nonlinear control, hybrid systems and stochastic optimization.
\end{IEEEbiography}
\end{document}